\documentclass[a4paper,11pt]{article}
\usepackage{braket,hyperref}
\usepackage{amssymb,mathtools,amsthm,amsmath}
\usepackage{tikz}
\usepackage{subfiles}
\usepackage{xargs}
\usepackage{enumitem}
\usepackage[colorinlistoftodos,prependcaption,textsize=tiny]{todonotes}

\theoremstyle{plain}
\newtheorem{theorem}{\bf Theorem}
\newtheorem{claim}[theorem]{\bf Claim}

\newtheorem{proposition}[theorem]{\bf Proposition}
\newtheorem{corollary}[theorem]{\bf Corollary}
\newtheorem{lemma}[theorem]{\bf Lemma}
\theoremstyle{definition}

\newenvironment{remark}[1][Remark.]{\begin{trivlist}
		\item[\hskip \labelsep {\bfseries #1}]}{\end{trivlist}}

\numberwithin{theorem}{section} 
\numberwithin{equation}{section}
\newcommandx{\error}[2][1=]{\todo[linecolor=red,backgroundcolor=red!25,bordercolor=red,#1]{#2}}
\newcommandx{\improvement}[2][1=]{\todo[linecolor=blue,backgroundcolor=blue!25,bordercolor=blue,#1]{#2}}
\newcommandx{\simania}[2][1=]{\todo[linecolor=green,backgroundcolor=green!25,bordercolor=green,#1]{#2}}

\newcommand{\Rea}{{\mathbb R}}
\newcommand{\QQ}{{\mathbb Q}}
\newcommand{\ZZ}{{\mathbb Z}}

\newcommand{\laplacian}[2]{L_{#1}\left(#2\right)}
\newcommand{\lappm}[3]{L_{#1}^{#3}\left(#2\right)}
\newcommand{\lap}[1]{L_{#1}}

\newcommand{\mineig}[2]{\mu_{#1}(#2)}

\newcommand{\homology}[2]{\tilde{H}_{#1}\left(#2;\Rea\right)}
\newcommand{\cohomology}[2]{\tilde{H}^{#1}\left(#2;\Rea\right)}

\DeclareMathOperator{\lk}{lk}
\newcommand{\degree}[1]{\deg_{#1}}

\newcommand{\normsquare}[1]{\left\| #1 \right\|^2}

\newcommand{\faces}[2]{#2(#1)}
\newcommand{\cochains}[2]{C^{#1}(#2)}
\newcommand{\sgn}[2]{(#1 : #2)}
\newcommand{\basiselement}[1]{1_{#1}}
\newcommand{\standardbasis}[1]{\{\basiselement{\sigma}\}_{\sigma\in #1}}
\newcommand{\orderedunion}[2]{[#1,#2]}

\DeclareMathOperator{\permsign}{sign}

\newcommand{\cobound}[2]{d_{#1}#2}
\newcommand{\bound}[2]{d_{#1}^*#2}
\newcommand{\matrixrep}[1]{\left[#1\right]}
\newcommand{\matrixrepel}[3]{\left[#1\right]_{#2,#3}}
\DeclareMathOperator{\Ker}{Ker}
\DeclareMathOperator{\Ima}{Im}
\newcommand{\inner}[2]{\left\langle#1,#2\right\rangle}
\newcommand{\norm}[1]{\left\|#1\right\|}
\newcommand{\missingfaces}{\mathcal{M}}
\newcommand{\missfaces}[2]{\mathcal{M}_{#2}(#1)}
\newcommand{\missingdims}[1]{J_{#1}}
\newcommand{\setsforrep}[1]{S(#1)}
\newcommand{\missingsubfaces}[1]{T(#1)}
\DeclareMathOperator{\mioperator}{Mis}
\newcommand{\mi}[1]{\mioperator\left(#1\right)}
\newcommand{\mis}[1]{m\left(#1\right)}
\newcommand{\cupdot}{\mathbin{\mathaccent\cdot\cup}}
\newcommand{\maxeig}[2]{\lambda_{\max}^{#1}(#2)}

\newcommand{\conn}[1]{\eta(#1)}
\DeclareMathOperator{\con}{conn}
\newcommand{\tot}[1]{\tilde{\gamma}(#1)}

\newcommand{\vrep}[1]{P_{#1}}

\newcommand{\tempmat}[2]{M_{#1,#2}}
\newcommand{\mata}[1]{A_{#1}}
\newcommand{\matb}[1]{B_{#1}}
\newcommand{\matc}[1]{R_{#1}}
\newcommand{\miscomplex}[1]{K_{#1}}
\newcommand{\maxmis}{h}
\newcommand{\eigc}[1]{\lambda_{k}}
\newcommand{\identity}{I}

\DeclareMathOperator{\matroidclosure}{cl}
\DeclareMathOperator{\matroidrank}{\rho}
\newcommand{\multiplier}[2]{#2_{#1}}
\newcommand{\varph}{\varphi_{M}}

\begin{document}
\title{Spectral gaps of simplicial complexes\\ without large missing faces}
\author{Alan Lew\footnote{Department of Mathematics, Technion, Haifa 32000, Israel. e-mail: alan@tx.technion.ac.il . Supported by ISF grant no. 326/16.}}

\date{}
\maketitle

\begin{abstract}
	
Let $X$ be a simplicial complex on $n$ vertices without missing faces of dimension larger than $d$. Let $\lap{j}$ denote the $j$-Laplacian acting on real $j$-cochains of $X$ and let $\mineig{j}{X}$ denote its minimal eigenvalue. We study the connection between the spectral gaps $\mineig{k}{X}$ for $k\geq d$ and $\mineig{d-1}{X}$. In particular, we establish the following vanishing result: If $\mineig{d-1}{X}>(1-\binom{k+1}{d}^{-1})n$, then $\cohomology{j}{X}=0$ for all $d-1\leq j \leq k$. As an application we prove a fractional extension of a Hall-type theorem of Holmsen, Mart\'inez-Sandoval and Montejano for general position sets in matroids.

%For a simplicial complex $X$ without missing faces of dimension larger than $d$ we study the connection between the spectral gaps of the $k$-Laplacian $\lap{k}$ for $k\geq d$ and the $(d-1)$-Laplacian $\lap{d-1}$. As an application we prove a fractional extension of a Hall type theorem of Holmsen, Martinez-Sandoval and Montejano for general position sets in matroids.
\end{abstract}

\section{Introduction}
\label{sec:intro}

Let $X$ be a simplicial complex on vertex set $V$. A simplex $\sigma\subset V$ is called a \emph{missing face} of $X$ if $\sigma\notin X$ but for any $\tau\subsetneq\sigma$, $\tau\in X$. 
The set of missing faces $\missingfaces_X$ of the complex $X$ completely determines $X$:
\[
X=\set{ \tau\subset V:\quad \sigma \not\subset\tau \text{ for all } \sigma\in\missingfaces_X}.
\]
Let $\maxmis(X)=\max \set{ \dim(\sigma) : \sigma\in \missingfaces_X}$.

For $k\geq -1$ let $\cochains{k}{X}$ be the space of real valued $k$-cochains of the complex $X$ and let $\cobound{k}:\cochains{k}{X}\to\cochains{k+1}{X}$ be the coboundary operator.
For $k\geq 0$ the reduced $k$-dimensional Laplacian of $X$ is defined by 
\[
\lap{k}(X)=\cobound{k-1}\bound{k-1}+\bound{k}\cobound{k}.
\]
$\lap{k}$ is a positive semidefinite operator from $\cochains{k}{X}$ to itself. The $k$-th \emph{spectral gap} of $X$, denoted by $\mineig{k}{X}$, is the smallest eigenvalue of $\lap{k}$.

Let $G=(V,E)$ be a graph on $n$ vertices. Its clique complex (or flag complex) $X(G)$ is the simplicial complex on vertex set $V$ whose simplices are the cliques of $G$.
Note that clique complexes are exactly the complexes with $\maxmis(X)=1$. Indeed, the missing faces of $X(G)$ are the edges of the complement of $G$. Aharoni, Berger and Meshulam \cite{aharoni2005eigenvalues} prove the following result:

\begin{theorem}[Aharoni, Berger, Meshulam \cite{aharoni2005eigenvalues}] %???
\label{thm:f1}
Let $G=(V,E)$ be a graph, where $|V|=n$, and let $X=X(G)$ be its clique complex. Then for $k\geq 1$
\[
	k \mineig{k}{X} \geq (k+1) \mineig{k-1}{X} - n.
\]
\end{theorem}

Our main result is a generalization of Theorem \ref{thm:f1} to complexes without large missing faces.
\begin{theorem}
\label{thm:fp}
Let $X$ be a simplicial complex with $\maxmis(X)=d$ on vertex set $V$, where $|V|=n$. Then for $k\geq d$
\[
	(k-d+1) \mineig{k}{X} \geq (k+1) \mineig{k-1}{X} - dn.
\]
\end{theorem}
Our proof combines the approach of \cite{aharoni2005eigenvalues} with additional new ideas. Both results can be thought of as global variants of Garland's method, which in its original form relates the spectral gaps of a complex with the spectral gaps of the links of its faces; See \cite{garland1973p,papikian2016garland}.  As a consequence of Theorem \ref{thm:fp} we obtain
\begin{theorem}
	\label{cor:fp}
	Let $X$ be a simplicial complex with $\maxmis(X)=d$, on vertex set $V$, where $|V|=n$. If
	\[
	\mineig{d-1}{X}>\left(1-\binom{k+1}{d}^{-1}\right)n,
	\]
	then
	$\cohomology{j}{X}=0$ for all $d-1\leq j\leq k$.
\end{theorem}

\begin{remark}
%For $d=1$ it is shown in \cite{ABM} that the assumption in Corollary \ref{cor:fp} cannot be replaced by $\mineig{0}{X}\geq (1-\frac{1}{k+1})n$: Let $n=r\ell$ for $r\geq 1, l\geq 2$, and let $G$ be the Tur\'an graph $T_r(n)$, i.e. the complete $r$-partite graph with all sides equal of size $\ell$. The clique complex $X=X(G)$ has $\mineig{0}{X}=\ell(r-1)=\frac{r-1}{r} n$, but $\cohomology{r-1}{X}\neq 0$.

%Currently we don't know if the assumption of the corollary is tight for $d\geq 1$, except for a few cases:

In the case $d=1$ it is shown in \cite{aharoni2005eigenvalues} that the condition in Theorem \ref{cor:fp} is the best possible: Let $G$ be the complete $r$-partite graph on $n=\ell r$ vertices, with all sides of size $\ell$. Then $\mineig{0}{X(G)}=\frac{r-1}{r} n$, but $\cohomology{r-1}{X(G)}\neq 0$.

For $d=2$ we have found such extremal examples only for a few cases:
\begin{enumerate}[leftmargin=*]
\item  Let $X$ be the simplicial complex whose vertices $V$ are the points of the affine plane over $\mathbb{F}_3$, and whose missing faces are the lines of the affine plane. On the one hand, one can check that $\mineig{1}{X}=6=\frac{2}{3}|V|$. On the other hand, $\cohomology{2}{X}=\Rea\neq 0$ (computer checked).

\item Let $X$ be the simplicial complex whose vertices $V$ are the points of the projective space of dimension $3$ over $\mathbb{F}_3$, and whose simplices are the sets of points containing at most two points from each line (so the missing faces are the subsets of size $3$ of the lines in the projective space). One can show that $\mineig{1}{X}=36=\frac{9}{10}|V|$. On the other hand, $\cohomology{4}{X}\neq 0$ (computer checked).
\end{enumerate}
\end{remark}

We next give some applications of Theorem \ref{thm:fp} to connectivity bounds and Hall-type theorems for general simplicial complexes.

Let $\conn{X}=\con_{\Rea}(X)+2$, where
\[
	\con_{\Rea}(X)=\min\set{i: \cohomology{i}{X}\neq 0}-1
\]
is the \emph{homological connectivity} of $X$ over $\Rea$.
%The \emph{homological connectivity} of $X$ over $\Rea$ is \[\con_{\Rea}(X)=\min\set{i: \cohomology{i}{X}\neq 0}-1.\]
%Let $\eta(X)=\con_{\Rea}(X)+2$.

%Let $G=(V,E)$ be a graph. 
A subset of vertices $S\subset V$ in a graph $G=(V,E)$ is called a \emph{totally dominating set} if for all $v\in V$ there is some $u\in S$ such that $vu\in E$. The \emph{total domination number} of $G$, denoted by $\tot{G}$, is the minimal size of a totally dominating set.
Let $I(G)$ be the independence complex of the graph, i.e. the simplicial complex whose faces are all the independent sets $\sigma\subset V$. The total domination number gives a lower bound on the connectivity of $I(G)$ (see \cite[Theorem 1.2]{meshulam2003domination}):         %(see \cite{meshulam2001clique}):
\begin{equation}
\label{eq:conn_vs_gamma_for_graphs}
\conn{I(G)}\geq \tot{G}/2.%\tag{$\ast$}
\end{equation}
(For additional lower bounds on $\conn{I(G)}$ in terms of other domination parameters, see e.g. \cite{aharoni2000hall,meshulam2003domination}).

The inequality \eqref{eq:conn_vs_gamma_for_graphs} had been generalized to general simplicial complexes: Let $X$ be a complex on vertex set $V$. We say that a subset $S\subset V$ is \emph{totally dominating} if for every $v\in V$ there is some $\sigma\subset S$  such that $\sigma\in X$ but $v\sigma\notin X$. The \emph{total domination number} of $X$, denoted $\tot{X}$, is the minimal size of a totally dominating set in $X$. For a graph $G$ we have $\tot{G}=\tot{I(G)}$ (the totally dominating sets of $I(G)$ are the same as the totally dominating sets of $G$).
In \cite{aharoni2006intersection} it is shown that for any simplicial complex $X$, $\conn{X}\geq \tot{X}/2$.

Another graphical domination parameter, $\Gamma(G)$, has been introduced in \cite{aharoni2005eigenvalues} as follows. A \emph{vector representation} of the graph $G$ is an assignment $P: V\to \Rea^{\ell}$ such that $P(v)\cdot P(w)\geq 1$ if $v$ and $w$ are adjacent in $G$, and $P(v)\cdot P(w)\geq 0$ otherwise. A non-negative vector $\alpha\in \Rea^V$ is called \emph{dominating} for $P$ if $\sum_{v\in V} \alpha(v) P(v) \cdot P(w) \geq 1$ for every $w\in V$. The \emph{value} of $P$ is
\[
|P|=\min\set{\sum_{v\in V} \alpha(v) : \,\alpha \text{ is dominating for $P$}}.
\]
Let $\Gamma(G)$ be the supremum of $|P|$ over all vector representations of $G$.
It is easy to check that $\Gamma(G)\leq \tilde{\gamma}(G)$ (see Proposition \ref{prop:gamma_vs_gamma}). In \cite{aharoni2005eigenvalues} the following was proved:
\begin{theorem}[Aharoni, Berger, Meshulam \cite{aharoni2005eigenvalues}]
\label{thm:gammavsetaABM}
\[
\conn{I(G)}\geq \Gamma(G).
\]	
\end{theorem}
With a view towards generalizing Theorem \ref{thm:gammavsetaABM} to an arbitrary simplicial complex $X$, we define a new domination parameter $\Gamma(X)$.

For $k\in\mathbb{N}$ let $\missfaces{k}{X}$ be the set of missing faces of $X$ of dimension $k$. Let
$\missingdims{X}=\set{ i\in \mathbb{N} : \, \missfaces{i}{X}\neq\emptyset}$
be the set of dimensions of simplices in $\missingfaces_X$. Define $\setsforrep{X}=\cup_{i\in \missingdims{X}} \binom{V}{i-1}$.

Let $\sigma\in\setsforrep{X}$ and fix $\ell=\ell(\sigma)\in \mathbb{N}$. A \emph{vector representation of $X$ with respect to $\sigma$} is an assignment $\vrep{\sigma}:V\to \Rea^{\ell}$
%$P_{\sigma}$ of a vector $\vrep{\sigma}(v)\in\Rea^{\ell}$ to each vertex $v\in V$,
such that the inner product $\vrep{\sigma}(v)\cdot \vrep{\sigma}(w)\geq 1$ if $vw\sigma\in \missfaces{|\sigma|+1}{X}$, and  $\vrep{\sigma}(v)\cdot \vrep{\sigma}(w)\geq 0$ otherwise. We identify the representation $\vrep{\sigma}$ with the matrix $\vrep{\sigma}\in\Rea^{|V|\times\ell}$ whose rows are the vectors $\vrep{\sigma}(v)$, for $v\in V$. We call the collection $P=\set{\vrep{\sigma}: \sigma\in\setsforrep{X}}$ a \emph{vector representation} of $X$.

For each $\sigma\in\setsforrep{X}$, let $\alpha_{\sigma}\in \Rea^{V}$ be a non-negative vector. The set $\set{\alpha_{\sigma} : \sigma\in\setsforrep{X}}$ is called \emph{dominating for $P$} if
\[
	\sum_{\sigma\in\setsforrep{X}} \alpha_{\sigma} \vrep{\sigma} \vrep{\sigma}^T \geq \textbf{1}
\] 
(where $\textbf{1}\in \Rea^{V}$ is the all $1$ vector). The \emph{value} of $P$ is
\[
	|P|=\min\left\{\sum_{\sigma\in\setsforrep{X}} \alpha_{\sigma}\cdot \textbf{1} :\, \set{\alpha_{\sigma}}_{\sigma\in\setsforrep{X}} \text{ is dominating for $P$}\right\}.
\]
Let $\Gamma(X)$ be the supremum of $|P|$ over all vector representations $P$ of $X$.

\begin{remark}[Remarks.]\leavevmode
\begin{enumerate}[leftmargin=*]
\item  If $X=I(G)$ for a graph $G$, then $\Gamma(X)$ coincides with the parameter $\Gamma(G)$ defined in \cite{aharoni2005eigenvalues}.

\item In the case when all the missing faces are of the same size, we can bound $\Gamma(X)$ by the total domination number $\tot{X}$:
\end{enumerate}
\begin{proposition}
\label{prop:gamma_vs_gamma}
Let $X$ be a simplicial complex with all its missing faces of dimension equal to $d$. Then
\[
	\Gamma(X)\leq \binom{\tot{X}}{d}.
\]
\end{proposition}
\end{remark}

Our main application of Theorem \ref{thm:fp} is the following extension of Theorem \ref{thm:gammavsetaABM}.
 
\begin{theorem}
\label{thm:connectivity_bound}
%Let $X$ be a simplicial complex on vertex set $V$, where $|V|=n$. Let $\missingfaces$ be the set of missing faces of $X$, and let $I=\set{i: \exists \sigma\in\missingfaces,\, \dim(\sigma)=i}$ be the set of dimensions of sets in $\missingfaces$. Then
%Let $X$ be a simplicial complex. Then
\[
\sum_{r\in \missingdims{X}} r \binom{\conn{X}}{r}\geq \Gamma(X).
\]
\end{theorem}
%In the case $X=I(G)$, then $I=\set{1}$, and we get $\conn{I(G)}\geq \Gamma(G)$ (this is Theorem 1.3 in \cite{ABM})

Let $V_1,\ldots,V_m$ be a partition of the vertex set $V$. We say that a subset $\sigma\subset V$ is \emph{colorful} if $|\sigma\cap V_i|=1$ for all $i\in\set{1,2,\ldots,m}$. Theorem \ref{thm:connectivity_bound} gives rise to the following Hall-type condition for the existence of colorful simplices:
\begin{theorem}
\label{thm:generalhalltype}
If for every $\emptyset\neq I\subset \set{1,2,\ldots,m}$
\[
	\Gamma(X[\cupdot_{i\in I} V_i])>
	\sum_{r\in \missingdims{X[\cupdot_{i\in I} V_i]}} r \binom{|I|-1}{r},
\]
then $X$ has a colorful simplex.
\end{theorem}

Next we show an application of Theorem \ref{thm:generalhalltype}.
Let $M$ be a matroid on vertex set $V$ with rank function $\matroidrank$. Assume $\rho(V)=d+1$. We identify $M$ with the simplicial complex of its independent sets.
For $S\subset V$, define its \emph{closure} by $\matroidclosure(S)=\set{v\in V: \, \matroidrank(S)=\matroidrank(S\cup\{v\})}.$
A subset $F\subset V$ is a \emph{flat} of $M$ if $F=\matroidclosure(F)$, i.e. $\matroidrank(F\cup\{v\})>\matroidrank(F)$ for all $v\notin F$.% A flat of rank $d$ is called a hyperplane. 

We say that a subset $S\subset V$ is in \emph{general position} with respect to $M$ if for any $1\leq k\leq d$ every flat of $M$ of rank $k$ contains at most $k$ points of $S$. This is equivalent to requiring that any $S'\subset S$ with $|S'|\leq d+1$ is an independent set in $M$.
% Note that if $S$ has size at least $d+1$, this is equivalent to requiring that any $d+1$ points in $S$ are independent in $M$, i.e. they don't all lie on a hyperplane.

For $S\subset V$ denote by $\varph(S)$ the maximal size of a subset of $S$ in general position.

Let $V_1,\ldots, V_m$ be a partition of $V$. The following Hall-type theorem is proved in \cite{holmsen2016geometric}.

\begin{theorem}[Holmsen, Mart\'inez--Sandoval, Montejano \cite{holmsen2016geometric} ]
\label{thm:HMS}
If for every $\emptyset\neq I\subset \set{1,2,\ldots,m}$
\[
	\varph(\cupdot_{i\in I} V_i)>
	\begin{cases}
		 |I|-1 & \text{ if } |I|\leq d+1,\\
		  d\binom{2|I|-2}{d} & \text{ if } |I|\geq d+2,
	\end{cases}
\]
then $V$ has a colorful subset in general position.
% with respect to $M$.
\end{theorem}

Let $S\subset V$. A weight function $f:S\to \Rea_{\geq 0}$ is in \emph{fractional general position} with respect to $M$ if for any $1\leq k\leq d$ and for any flat $F$ of $M$ of rank $k$ and $\sigma\subset F\cap S$ of size $k-1$,
\[
	\sum_{\substack{v\in S,\\ \matroidclosure(v\sigma)=F}} f(v) \leq d.
\]
	
Denote by $\varph^*(S)$ the maximum of $\sum_{v\in S}f(v)$ over all functions $f:S\to \Rea_{\geq 0}$ in fractional general position.
Let $f$ be the characteristic function of a set $S'\subset S$ in general position. Let $F$ be a flat of $M$ of rank $k$ for $1\leq k\leq d$ and $\sigma\subset F\cap S$ of size $k-1$. Then
\[
	\sum_{\substack{v\in S,\\ \matroidclosure(v\sigma)=F}} f(v)
	=\left| \set{ v\in S' : \, \matroidclosure(v\sigma)=F} \right|\leq |S'\cap F| \leq k\leq d,
\]
so $f$ is in fractional general position. Therefore
\begin{equation}
\label{eq:fractionalphi}
\varph^*(S)\geq \varph(S).
\end{equation}
Here we prove the following:% strengthening of Theorem \ref{thm:HMS}:

\begin{theorem}
	\label{thm:myHMSstar}
	If for every $\emptyset\neq I\subset \set{1,2,\ldots,m}$
	\[
	\varph^*(\cupdot_{i\in I} V_i)>
	d \sum_{r=1}^d r\binom{|I|-1}{r},
	\]
	then $V$ contains a colorful subset in general position.
\end{theorem}

In particular, we obtain a strengthening of Theorem \ref{thm:HMS}:
\begin{theorem}
\label{thm:myHMS}
	If for every $\emptyset\neq I\subset \set{1,2,\ldots,m}$
	\[
	\varph(\cupdot_{i\in I} V_i)>
		\begin{cases}
		|I|-1 & \text{ if } |I|\leq d+1,\\
		 d \sum_{r=1}^d r\binom{|I|-1}{r} & \text{ if } |I|\geq d+2,
		\end{cases}
	\]
	then $V$ contains a colorful subset in general position.
\end{theorem}

The paper is organized as follows. In Section \ref{sec:prel} we review some basic facts concerning simplicial cohomology and high dimensional Laplacians. We also introduce some notation and results about complexes without large missing faces that we will need later. In Section \ref{sec:spectral} we prove our main result, Theorem \ref{thm:fp}, and its corollary Theorem \ref{cor:fp}. Section \ref{sec:vector_domination} deals with the vector domination parameter $\Gamma(X)$ of the complex $X$. In it we prove Proposition \ref{prop:gamma_vs_gamma}, Theorem \ref{thm:connectivity_bound} and Theorem \ref{thm:generalhalltype}. In Section \ref{sec:matroids} we apply the results of the previous section in order to prove Theorems \ref{thm:myHMSstar} and \ref{thm:myHMS}, which provide sufficient conditions for the existence of colorful sets in general position in a matroid.

\section{Preliminaries}
\label{sec:prel}

\subsection{Simplicial cohomology}
Let $X$ be a finite simplicial complex on the vertex set $V$. We denote the set of $k$-dimensional simplices in $X$ by $\faces{k}{X}$. For each $\sigma\in X(k)$ we choose an order of its vertices $v_0,\ldots,v_k$, which induces an orientation on $\sigma$.

For $\sigma\in \faces{k}{X}$, let
$
	\lk(X,\sigma)=\set{ \tau\in X: \tau\cup \sigma\in X, \tau\cap\sigma=\emptyset}
$
be the \emph{link} of $\sigma$ in $X$, and 
$\degree{X}(\sigma)=|\set{\eta\in\faces{k+1}{X}: \sigma\subset \eta}|$ be the \emph{degree} of $\sigma$ in $X$.
For $U\subset V$, let $X[U]=\set{\sigma\in X: \sigma\subset U}$ be the subcomplex of $X$ induced by $U$.

For two ordered simplices $\sigma\in X$, $\tau \in \lk(X,\sigma)$, denote by $\orderedunion{\sigma}{\tau}$, or simply by $\sigma\tau$,  their ordered union. Similarly, for $v\in V$ denote by $v\sigma$ the ordered union of $\{v\}$ and $\sigma$. %If $v\in V$ (not necessarily in the link of $\sigma$) we denote by $v\sigma$ the union $\set{v}\cup \sigma$ (without assigning an order of the vertices).
%\improvement{explain better this abuse of notation.}

For $\sigma\in X$, and $\tau\subset \sigma$, both given an order on their vertices, we define $\sgn{\sigma}{\tau}$ to be the sign of the permutation on the vertices of $\sigma$ which maps the ordered simplex $\sigma$ to the ordered simplex $\orderedunion{\sigma\setminus\tau}{\tau}$ (where the order on the vertices of $\sigma\setminus\tau$ is the one induced by the order on $\sigma$).
%\improvement[inline]{change 'ordering' to just 'order'?}
%For $\tau=\sigma$ (as sets), we say that $\sigma$ and $\tau$ have the same orientation if $\sgn{\sigma}{\tau}=1$ , and opposite orientation if $\sgn{\sigma}{\tau}=-1$.

A simplicial $k$-cochain is a real valued skew-symmetric function on all ordered $k$-simplices. That is, $\phi$ is a $k$-cochain if for any two $k$-simplices $\sigma,\tilde{\sigma}$ in $X$ that are equal as sets, it satisfies $\phi(\tilde{\sigma})=\sgn{\tilde{\sigma}}{\sigma} \phi(\sigma)$.

For $k\geq 0$ let $\cochains{k}{X}$ denote the space of $k$-cochains on $X$. For $k=-1$ we define $\cochains{-1}{X}=\Rea$.

We will use the following lemma implicitly in future calculations.
\begin{lemma}
%Let $\tau,\eta \subseteq\sigma$ and $\theta\subseteq\tau\cap\eta$ be ordered simplices, and let $\tilde{\sigma}$, $\tilde{\tau}$, $\tilde{\eta}$, $\tilde{\theta}$ be equal as sets to $\sigma$, $\tau$, $\eta$ and $\theta$ respectively. Assume $\eta,\tau\in X(k)$, and let $\phi$ be a $k$-cochain on $X$. Then
Let $\tau,\eta\in\faces{k}{X}$ and $\phi\in\cochains{k}{X}$. Let $\sigma,\theta\in X$ be ordered simplices such that $\tau,\eta\subset \sigma$ and $\theta\subset \tau\cap\eta$, and let $\tilde{\sigma}$, $\tilde{\tau}$, $\tilde{\eta}$, $\tilde{\theta}$ be equal as sets to $\sigma$, $\tau$, $\eta$ and $\theta$ respectively. Then

\begin{enumerate}[leftmargin=*]
 \item %Let $\tau\subset\sigma$ be ordered simplices, and let $\tilde{\sigma}$, $\tilde{\tau}$ be equal as sets to $\sigma$ and $\tau$ respectively. Then
$
     \sgn{\sigma}{\tau}=\sgn{\sigma}{\tilde{\tau}}\cdot \sgn{\tilde{\tau}}{\tau},
$
and if $|\sigma\setminus\tau|=1$ then
$
    \sgn{\sigma}{\tau}=\sgn{\sigma}{\tilde{\sigma}}\cdot \sgn{\tilde{\sigma}}{\tau}.
$
\item  $\phi(\tau)^2=\phi(\tilde{\tau})^2.$
\item $\sgn{\sigma}{\tau}\phi(\tau)=\sgn{\sigma}{\tilde{\tau}}\phi(\tilde{\tau}), $ and if $|\tau\setminus\theta|=1$ then \[\sgn{\tau}{\theta}\phi(\tau)=\sgn{\tilde{\tau}}{\theta}\phi(\tilde{\tau}).\]
\item If $|\sigma\setminus\tau|=1$ and $|\sigma\setminus\eta|=1$ then
\[\sgn{\sigma}{\tau}\sgn{\sigma}{\eta}\phi(\tau)\phi(\eta)=\sgn{\tilde{\sigma}}{\tilde{\tau}}\sgn{\tilde{\sigma}}
{\tilde{\eta}}\phi(\tilde{\tau})\phi(\tilde{\eta}).\]
\item If $|\tau\setminus\theta|=1$ and $|\eta\setminus\theta|=1$ then
\[\sgn{\tau}{\theta}\sgn{\eta}{\theta}\phi(\tau)\phi(\eta)=\sgn{\tilde{\tau}}{\tilde{\theta}}\sgn{\tilde{\eta}}
{\tilde{\theta}}\phi(\tilde{\tau})\phi(\tilde{\eta}).\]
\end{enumerate}
\end{lemma}
\begin{proof}\
\begin{enumerate}[leftmargin=*]
\item
Let $\pi_1$ be the permutation on the vertices of $\sigma$ that maps $\sigma$ to $\orderedunion{\sigma\setminus\tilde{\tau}}{\tilde{\tau}}=\orderedunion{\sigma\setminus\tau}{\tilde{\tau}}$, and let $\pi_2$ be the permutation on the vertices of $\tau$ that maps $\tilde{\tau}$ to $\tau$. Extend $\pi_2$ to a permutation $\tilde{\pi}_2$ on the vertices of $\sigma$, which maps $\orderedunion{\sigma\setminus\tau}{\tilde{\tau}}$ to $\orderedunion{\sigma\setminus\tau}{\tau}$. It satisfies $\permsign(\pi_2)=\permsign(\tilde{\pi}_2)$. Define $\pi=\tilde{\pi}_2\circ \pi_1$. $\pi$ maps $\sigma$ to $\orderedunion{\sigma\setminus\tau}{\tau}$, therefore
\begin{multline*}
\sgn{\sigma}{\tau}=\permsign(\pi)
%%%%%%%%%%%%%
=\permsign(\tilde{\pi}_2)\cdot\permsign(\pi_1)\\
%%%%%%%%%%%%%
=\permsign(\pi_2)\cdot\permsign(\pi_1)
%%%%%%%%%%%%
=\sgn{\tilde{\tau}}{\tau}\cdot\sgn{\sigma}{\tilde{\tau}}.
\end{multline*}
Assume now that $|\sigma\setminus\tau|=1$ and let $\{v\}=\sigma\setminus\tau$. Let $\pi_3$ be the permutation on the vertices of $\sigma$ that maps $\sigma$ to $\tilde{\sigma}$, and $\pi_4$ be the permutation which maps $\tilde{\sigma}$ to $\orderedunion{\tilde{\sigma}\setminus\tau}{\tau}= v\tau= \orderedunion{\sigma\setminus\tau}{\tau}$. Then the permutation $\pi'=\pi_4 \circ \pi_3$ maps $\sigma$ to $\orderedunion{\sigma\setminus\tau}{\tau}$, therefore
\[
\sgn{\sigma}{\tau}=\permsign(\pi')=\permsign(\pi_4)\cdot\permsign(\pi_3)=
\sgn{\sigma}{\tilde{\sigma}}\cdot\sgn{\tilde{\sigma}}{\tau}.
\]

\item Since $\phi$ is a cochain, we have $\phi(\tau)^2=\sgn{\tau}{\tilde{\tau}}^2 \phi(\tilde{\tau})^2=\phi(\tilde{\tau})^2.$
\item By the first part of this lemma
\[
\sgn{\sigma}{\tau}\phi(\tau)=\sgn{\sigma}{\tilde{\tau}}
\sgn{\tilde{\tau}}{\tau}\phi(\tau),
\]
and since $\phi$ is a cochain
\[
\sgn{\sigma}{\tilde{\tau}}
\sgn{\tilde{\tau}}{\tau}\phi(\tau)=
 \sgn{\sigma}{\tilde{\tau}}\phi(\tilde{\tau}).
\]

The second equality is similar: By the first part of the lemma
\[
\sgn{\tau}{\theta}\phi(\tau)=\sgn{\tau}{\tilde{\tau}}
\sgn{\tilde{\tau}}{\theta}\phi(\tau),
\]
and since $\phi$ is a cochain
\[
\sgn{\tau}{\tilde{\tau}}
\sgn{\tilde{\tau}}{\theta}\phi(\tau)=\sgn{\tau}{\tilde{\tau}}
\sgn{\tilde{\tau}}{\theta}\sgn{\tau}{\tilde{\tau}}\phi(\tilde{\tau})
=\sgn{\tilde{\tau}}{\theta}\phi(\tilde{\tau}).
\]

\item By part 3 of this lemma we have
\[
\sgn{\sigma}{\tau}\sgn{\sigma}{\eta}\phi(\tau)\phi(\eta)=
\sgn{\sigma}{\tilde{\tau}}\sgn{\sigma}{\tilde{\eta}}\phi(\tilde{\tau})\phi(\tilde{\eta}).
\]
Then by part 1
\begin{multline*}
\sgn{\sigma}{\tilde{\tau}}\sgn{\sigma}{\tilde{\eta}}\phi(\tilde{\tau})\phi(\tilde{\eta})= \sgn{\sigma}{\tilde{\sigma}}\sgn{\tilde{\sigma}}{\tilde{\tau}}\sgn{\sigma}{\tilde{\sigma}}\sgn{\tilde{\sigma}}{\tilde{\eta}}\phi(\tilde{\tau})\phi(\tilde{\eta})=\\
\sgn{\tilde{\sigma}}{\tilde{\tau}}\sgn{\tilde{\sigma}}
{\tilde{\eta}}\phi(\tilde{\tau})\phi(\tilde{\eta}).
\end{multline*}
\item The proof is similar to the proof of part 4.
\end{enumerate}
\end{proof}

For $k\geq 0$ let the coboundary operator $\cobound{k}{}: \cochains{k}{X}\to\cochains{k+1}{X}$ be the linear operator defined by
\[
    \cobound{k}{\phi}(\sigma)=\sum_{i=0}^{k+1} (-1)^i \phi(\sigma_i),
\]
where for an ordered $(k+1)$-simplex $\sigma=[v_0,\ldots, v_{k+1}]$, $\sigma_i$ is the ordered simplex obtained by removing the vertex $v_i$, that is $\sigma_i=[v_0,\ldots,\hat{v_i},\ldots, v_{k+1}]$.
Equivalently, we can write
%\begin{equation}
\[
%\label{eq:coboundary}
    \cobound{k}{\phi}(\sigma)=\sum_{\tau\in\sigma(k)}\sgn{\sigma}{\tau} \phi(\tau),
\]
%\end{equation}
where $\sigma(k)\subset \faces{k}{X}$ is the set of all $k$-dimensional faces of $\sigma$, each given some fixed order on its vertices.

For $k=-1$ we define $\cobound{-1}{}: \cochains{-1}{X}=\Rea \to\cochains{0}{X}$ by $\cobound{-1}{a}(v)=a$, for every $a\in \Rea$, $v\in V$.

Let $\cohomology{k}{X}=\Ker(\cobound{k}{})/\Ima(\cobound{k-1}{})$ be the $k$-th reduced cohomology  group of $X$ with real coefficients. %This is a vector space over $\Rea$, and we denote $\betti{k}{X}=\dim \cohomology{k}{X}$.

\subsection{Higher Laplacians}
For each $k\geq -1$ we define an inner product on $\cochains{k}{X}$ by
\[
\inner{\phi}{\psi}=\sum_{\sigma\in \faces{k}{X}} \phi(\sigma)\psi(\sigma).
\]
This induces a norm on $\cochains{k}{X}$:
\[
\norm{\phi}=\left(\sum_{\sigma\in\faces{k}{X}} \phi(\sigma)^2\right)^{1/2}.
\]

Let $\bound{k}{}: \cochains{k+1}{X}\to\cochains{k}{X}$ be the adjoint of $\cobound{k}{}$ with respect to this inner product. We can write $\bound{k}{\phi}$ explicitly:
%\begin{equation}
%\label{eq:boundary}
\[
    \bound{k}{\phi}(\tau)=\sum_{v\in \lk(X,\tau)}\phi(v\tau).
\]
%\end{equation}
For $k\geq 0$ define the lower $k$-Laplacian of $X$ by $\lappm{k}{X}{-}=\cobound{k-1}{}\bound{k-1}{}$ and the upper $k$-Laplacian of $X$ by 
$
\lappm{k}{X}{+}=\bound{k}{}\cobound{k}{}.
$
The reduced $k$-Laplacian of $X$ is the positive semidefinite operator on $\cochains{k}{X}$ given by $\laplacian{k}{X}=\lappm{k}{X}{-}+\lappm{k}{X}{+}$.

Let $k\geq 0$ and $\sigma\in \faces{k}{X}$. We define the $k$-cochain $\basiselement{\sigma}$ by
\[
    \basiselement{\sigma}(\tau)=\begin{cases}
    \sgn{\sigma}{\tau} & \text{ if $\sigma=\tau$ (as sets)},\\
    0 & \text{ otherwise.}
    \end{cases}
\]
The set $\standardbasis{\faces{k}{X}}$ forms a basis of the space $\cochains{k}{X}$, which we will call the standard basis.

For a linear operator $T: \cochains{k}{X}\to \cochains{k}{X}$, let $\matrixrep{T}$ be the matrix representation of $T$ with respect to the standard basis.
%basis $\standardbasis{\faces{k}{X}}$.
We denote by $\matrixrepel{T}{\sigma}{\tau}$ the matrix element of $\matrixrep{T}$ at index $(\basiselement{\sigma},\basiselement{\tau})$.

One can write explicitly the matrix representation of the Laplacian operators in the standard basis (see e.g. \cite{duval2002shifted,goldberg2002combinatorial}):
\begin{claim}%[Matrix representation of Laplacian operators]
\label{claim:lapmatrix}

For $k\geq 0$
\begin{align*}
&\matrixrepel{\lap{k}^{-}}{\sigma}{\tau}=
\begin{cases}
	k+1 & \mbox{if } \sigma = \tau, \\
	\sgn{\sigma}{\sigma\cap\tau}\cdot\sgn{\tau}{\sigma\cap\tau}
	& \mbox{if } \left|\sigma \cap \tau \right|=k,\\
	0 & \mbox{otherwise,}
\end{cases}	
\\
%\]
%\[
&\matrixrepel{\lap{k}^{+}}{\sigma}{\tau}=
\begin{cases}
	\degree{X}(\sigma) & \mbox{if } \sigma = \tau,\\
	-\sgn{\sigma}{\sigma \cap \tau}\cdot\sgn{\tau}{\sigma \cap \tau}
	& \mbox{if } \left|\sigma \cap \tau \right|=k,  \sigma \cup \tau \in X(k+1), \\
	0 & \mbox{otherwise,}
\end{cases}	
\end{align*}
and
\[
\matrixrepel{\lap{k}}{\sigma}{\tau}=
\begin{cases}
	k+1+\degree{X}(\sigma) & \mbox{if } \sigma = \tau, \\
	\sgn{\sigma}{\sigma \cap \tau}\cdot\sgn{\tau}{\sigma \cap \tau}
	& \mbox{if } \left|\sigma \cap \tau \right|=k,  \sigma \cup \tau \notin X(k+1), \\
	0 & \mbox{otherwise.}
\end{cases}	
\]	
\end{claim}

The following upper bound on the eigenvalues of the Laplacian is implicit in \cite{duval2002shifted}:
\begin{lemma}
	\label{lemma:max_eig}
	Let $X$ be a simplicial complex on vertex set $V$, with $|V|=n$. Let $k\geq 0$ and let $\lambda$ be an eigenvalue of $\laplacian{k}{X}$. Then 
	\[
	\lambda\leq n.
	\]
\end{lemma}

The following discrete version of Hodge's theorem had been observed by Eckmann in \cite{eckmann1944harmonische}.

\begin{theorem}[Simplicial Hodge theorem]
\[
	\cohomology{k}{X}\cong \Ker \lap{k}.
\]
\end{theorem}

As a consequence of Hodge theorem we obtain
\begin{corollary}
\label{cor:hodge}
$\cohomology{k}{X}=0$ if and only if $\mineig{k}{X}>0$.
\end{corollary}
%\begin{proof}
%By Hodge's theorem,  $\cohomology{k}{X}=0$ if and only if $\Ker \laplacian{k}{X}=0$. This happens if and only if $0$ is not an eigenvalue of $\laplacian{k}{X}$, which is the same as saying that all the eigenvalues of $\laplacian{k}{X}$  are positive, or equivalently, that $\mineig{k}{X}>0$.
%\end{proof}

\subsection{Missing faces and sums of degrees}

Let $X$ be a complex on vertex set $V$ with $\maxmis(X)=d$. Let $k\geq d$ and $\theta\in \binom{V}{k+1}$. Define

\[
\missingsubfaces{\theta}= \left\{ \tau\in \binom{\theta}{d+1}  : \tau\notin X(d)\right\}.
\]
So $\missingsubfaces{\theta}$ is the set of all $d$-dimensional simplices in $\theta$ that do not belong to $X$, and 
 $\theta\in X$ if and only if $\missingsubfaces{\theta}=\emptyset$.
Let
\[
\mi{\theta}= \bigcap_{\tau\in \missingsubfaces{\theta}} \tau
\]
and
\[
\mis{\theta}= \left|\bigcap_{\tau\in T(\theta)} \tau \right|.
\]
Since every $\tau\in T(\theta)$ has $d+1$ vertices it follows that $\mis{\theta}\leq d+1$. Another simple observation is the following:

\begin{lemma}
\label{lem:missing_faces}
	Let $\sigma,\tau\in \faces{k}{X}$ such that $|\tau\cap\sigma|=k$. Then if $\sigma\cup\tau\in \faces{k+1}{X}$, $\mis{\sigma\cup\tau}=0$, otherwise $2\leq\mis{\sigma\cup\tau}\leq d+1$.
\end{lemma}
\begin{proof}
Denote $\sigma\setminus \tau=\{v\}$ and $\tau\setminus\sigma=\{w\}$. If $\sigma\cup\tau\in X(k+1)$ then $T(\sigma\cup\tau)=\emptyset$, therefore $\mis{\sigma\cup\tau}=0$. If  $\sigma\cup\tau\notin X(k+1)$, then every $\eta\in T(\sigma\cup\tau)$ must contain both $v$ and $w$ (otherwise $\eta$ will be contained in $\sigma$ or in $\tau$, a contradiction to $\eta\notin 
\faces{d}{X}$). Therefore $\mis{\sigma\cup\tau}\geq 2$.	
\end{proof}

The following is a known result about clique complexes (see \cite[Claim 3.4]{aharoni2005eigenvalues}, \cite{aigner2001turan}):
\begin{lemma}
Let $X$ be a clique complex with $n$ vertices and let $\sigma\in X(k)$. Then
\[
\sum_{\tau\in\faces{k-1}{\sigma}} \degree{X}(\tau) - k\degree{X}(\sigma)\leq n.
\]
\end{lemma}

We will need a version of this lemma for complexes without large missing faces:
\begin{lemma}
\label{lemma:countdegrees}
Let $X$ be a simplicial complex on vertex set $V$ with $\maxmis(X)=d$. Let $k\geq d$ and $\sigma\in X(k)$. Then
\begin{align*}
	\sum_{\tau\in\faces{k-1}{\sigma}} \degree{X}(\tau)
	&=k+1+(k+1)\degree{X}(\sigma) \\
	&+\sum_{r=2}^{d+1} (r-1)\cdot\left|\set{v\in V: \mis{v\sigma}=r}\right|.
\end{align*}
\end{lemma}

\begin{proof}
\begin{multline}
\label{eq:degsum}
\sum_{\tau\in\faces{k-1}{\sigma}} \degree{X}(\tau)= \sum_{\tau\in \sigma(k-1)} \sum_{v\in \lk(X,\tau)} 1
%%%%%%%%
=\sum_{v\in V} \sum_{\substack{\tau\in \sigma(k-1),\\ \tau\in \lk(X,v)}}1\\
%%%%%%%%
=\sum_{v\in \sigma}\sum_{\substack{\tau\in \sigma(k-1),\\ \tau\in \lk(X,v)}}1+
\sum_{v\in\lk(X,\sigma)} \sum_{\substack{\tau\in \sigma(k-1),\\ \tau\in \lk(X,v)}}1+
\sum_{\substack{v\in V\setminus\sigma,\\ v\notin\lk(X,\sigma)}}\sum_{\substack{\tau\in \sigma(k-1),\\ \tau\in \lk(X,v)}}1.
\end{multline}
We consider separately the three summands on the right hand side of \eqref{eq:degsum}:
\begin{enumerate}[leftmargin=*]
\item For $v\in\sigma$, there is only one $\tau\in\faces{k-1}{\sigma}$ such that $\tau\in\lk(X,v)$, namely $\tau=\sigma\setminus\{v\}$. Thus the first summand is $k+1$.

\item For $v\in\lk(X,\sigma)$, any $\tau\in\faces{k-1}{\sigma}$ is in $\lk(X,v)$, therefore the second summand is $(k+1)\degree{X}(\sigma)$. 

\item Let $v\in V\setminus\sigma$ such that $v\notin\lk(X,\sigma)$. Let $\tau\in\faces{k-1}{\sigma}$ and let $u$ be the unique vertex in $\sigma\setminus\tau$. If $\tau\in\lk(X,v)$ then every missing face of $X$ contained in $v\sigma$ must contain $u$, so $u\in\mi{v\sigma}$. If $\tau\notin \lk(X,v)$, then there is a missing face of $X$ contained in $v\tau$, and therefore it doesn't contain the vertex $u$. Hence, $u\notin \mi{v\sigma}$. 
%For $v\in V\setminus(\sigma\cup\lk(X,\sigma))$, $\tau\in\faces{k-1}{\sigma}$ is in $\lk(X,v)$ if and only if the vertex $u=\sigma\setminus\tau$ is in $\mi{v\sigma}$.
Since $v\in\mi{v\sigma}$, the number of $\tau\in\faces{k-1}{\sigma}$ such that $\tau\in\lk(X,v)$ is exactly  $\mis{v\sigma}-1$.
Hence the third summand is 
\[
\sum_{\substack{v\in V\setminus\sigma,\\ v\notin\lk(X,\sigma)}}(\mis{v\sigma}-1) =\sum_{r=2}^{d+1} (r-1)|\set{v\in V :\, \mis{v\sigma}=r}|.
\]
\end{enumerate}
We obtain
\begin{multline*}
\sum_{\tau\in\faces{k-1}{\sigma}} \degree{X}(\tau)
=k+1+(k+1)\degree{X}(\sigma)\\+\sum_{r=2}^{d+1} (r-1)\left|\set{v\in V: \mis{v\sigma}=r}\right|.
\end{multline*}

\end{proof}

\section{Spectral gaps}
\label{sec:spectral}

In this section we prove Theorems \ref{thm:fp} and \ref{cor:fp}.

Let $X$ be a simplicial complex with $\maxmis(X)=d$ on vertex set $V$, where $|V|=n$, and let $k\geq d$. For $\phi\in\cochains{k}{X}$ and $u\in V$ we define $\phi_u\in\cochains{k-1}{X}$ by
\[
\phi_u(\tau)=
\begin{cases}
\phi(u\tau) & \mbox{if } u\in \lk(X,\tau),\\
0 & \mbox{otherwise.}
\end{cases}
\]
Let $\matb{k}:C^k(X)\to C^k(X)$ be the linear transformation whose matrix representation in the standard basis is% basis $\set{\varphi_\sigma}_{\sigma\in X(k)}$ is
\[
\matrixrepel{\matb{k}}{\tau}{\sigma}=
\begin{cases}
k \degree{X}(\sigma)-\sum_{\eta\in\faces{k-1}{\sigma}} \degree{X}(\eta) & \mbox{if } \sigma=\tau,\\
(\mis{\sigma\cup\tau}-2)\sgn{\sigma}{\sigma\cap\tau}\cdot\sgn{\tau}{\sigma\cap\tau} & \mbox{if } \substack{|\sigma\cap\tau|=k,\\ \sigma\cup\tau\notin X(k+1),}\\
0 & \mbox{otherwise.}
\end{cases}
\]
Let $\matc{k}= (d-1)\lap{k}-\matb{k}$, and let $\eigc{k}$ be the largest eigenvalue of $\matc{k}$.

The proof of Theorem \ref{thm:fp} depends on the following two ingredients:

\begin{proposition}
\label{prop:laplacian_identity}
Let $\phi\in\cochains{k}{X}$. Then
\[
(k-d+1)\inner{\lap{k}\phi}{\phi}= \sum_{u\in V} \inner{L_{k-1}\phi_u}{\phi_u}-\inner{\matc{k}\phi}{\phi}.
\]
\end{proposition}

\begin{proposition}
\label{prop:eigenbound}
$\eigc{k}\leq dn$.
\end{proposition}
We postpone the proof of these propositions to the end of this section, and first show how they imply Theorem \ref{thm:fp}.

\begin{proof}[Proof of Theorem \ref{thm:fp}]
Let $0\neq\phi\in \cochains{k}{X}$ be an eigenvector of $\lap{k}$ with eigenvalue $\mineig{k}{X}$. By Proposition \ref{prop:laplacian_identity} we obtain
\begin{multline*}
(k-d+1)\mineig{k}{X} \norm{\phi}^2
%%%%%
=(k-d+1)\inner{\lap{k}\phi}{\phi}\\
%%%%%%
=\sum_{u\in V} \inner{L_{k-1}\phi_u}{\phi_u}-\inner{\matc{k}\phi}{\phi}
%%%%%%
\geq \mineig{k-1}{X}\sum_{u\in V}\norm{\phi_u}^2 -\eigc{k}\norm{\phi}^2.
\end{multline*}
But
\begin{multline*}
\sum_{u\in V} \norm{\phi_u}^2
%%%%%
=\sum_{u\in V} \sum_{\tau\in \faces{k-1}{X}} \phi_u(\tau)^2
%%%%%
=\sum_{\tau\in\faces{k-1}{X}} \sum_{u\in\lk(X,\tau)} \phi(u\tau)^2\\
%%%%%
=(k+1)\sum_{\sigma\in \faces{k}{X}} \phi(\sigma)^2
%%%%%
=(k+1)\norm{\phi}^2.
\end{multline*}
Therefore
\[
(k-d+1)\mineig{k}{X} 
%%%%%
\geq (k+1) \mineig{k-1}{X} -\eigc{k},
\]
and by Proposition \ref{prop:eigenbound}
\[
(k-d+1)\mineig{k}{X} 
%%%%%
\geq (k+1) \mineig{k-1}{X} -dn.
\]
\end{proof}

For the proof of Theorem \ref{cor:fp} we will need the following result, which will also be used in Section \ref{sec:vector_domination}.
\begin{claim}
\label{claim:induction}
For $k\geq d-1$,
\begin{equation}
\label{eq:induction}
\mineig{k}{X} \geq \binom{k+1}{d} \mineig{d-1}{X}- \left( \binom{k+1}{d}-1\right)n.
\end{equation}
If in addition $X$ has complete $(d-1)$-dimensional skeleton, then there is equality in \eqref{eq:induction} for $0\leq k\leq d-1$.
%the inequality holds also for $0\leq k\leq d$ (and it actually achieves equality for all $0\leq k\leq d$).
\end{claim}
\begin{proof}
We argue by induction on $k$. The case $k=d-1$ is clear. Let $k\geq d$. By Theorem \ref{thm:fp} and the induction hypothesis we obtain
\begin{multline*}
\mineig{k}{X}\geq \frac{k+1}{k-d+1} \mineig{k-1}{X} - \frac{d}{k-d+1} n\\
%%%%%%%%%%%%
\geq  \frac{k+1}{k-d+1} \left[ \binom{k}{d} \mineig{d-1}{X}- \left( \binom{k}{d}-1\right)n\right] - \frac{d}{k-d+1} n \\
%%%%%%%%%%%%%
%= \binom{k+1}{d} \mineig{d-1}{X} - \binom{k+1}{d} n +\frac{k+1}{k-d+1}n-  \frac{d}{k-d+1} n\\
%%%%%%%%%%%%%%%%
= \binom{k+1}{d} \mineig{d-1}{X}- \left( \binom{k+1}{d}-1\right)n.
\end{multline*}

Now assume that $X$ has complete $(d-1)$-dimensional skeleton, and let $k<d-1$. Then we have $\binom{k+1}{d}=0$, therefore the inequality in the claim is just $\mineig{k}{X}\geq n$. But one can see by Claim \ref{claim:lapmatrix} that in this case $\lap{k}$ is the scalar matrix with diagonal elements $n$, thus $\mineig{k}{X}=n$.
\end{proof}

\begin{proof}[Proof of Theorem \ref{cor:fp}]
Let $d-1\leq j\leq k$. We have by Claim \ref{claim:induction}
\begin{align*}
\mineig{j}{X}&\geq \binom{j+1}{d}\mineig{d-1}{X}- \left( \binom{j+1}{d}-1\right)n\\
%%%%%%%%%%	
&>\binom{j+1}{d}\cdot\left(1- \binom{k+1}{d}^{-1}\right)n- \left( \binom{j+1}{d}-1\right)n\\
%%%%%%%%%%
&\geq \binom{j+1}{d}\cdot\left(1- \binom{j+1}{d}^{-1}\right)n- \left( \binom{j+1}{d}-1\right)n =0.
\end{align*}
Thus, by Corollary \ref{cor:hodge}, $\cohomology{j}{X}=0$.
\end{proof}

%Next we establish Proposition \ref{prop:laplacian_identity}:
In order to prove Proposition \ref{prop:laplacian_identity} we will need the following claims.

\begin{claim}[{see \cite[Claim 3.1]{aharoni2005eigenvalues}}]
\label{claim:dk}
For $\phi\in \cochains{k}{X}$
\[
	\norm{\cobound{k}{\phi}}^2= \sum_{\sigma\in \faces{k}{X}} \degree{X}(\sigma) \phi(\sigma)^2
	-2 \sum_{\eta\in \faces{k-1}{X}} \sum_{vw\in \lk(X,\eta)} \phi(v\eta) \phi (w\eta).
\]
\end{claim}
\begin{proof}
\begin{align*}
	\norm{\cobound{k}{\phi}}^2&
%%%%%
	=\sum_{\tau\in \faces{k+1}{X}} \cobound{k}{\phi(\tau)}^2
	\\
%%%%%
	&=\sum_{\tau\in \faces{k+1}{X}} \left(\sum_{\theta_1\in\faces{k}{\tau}} \sgn{\tau}{\theta_1} \phi(\theta_1)\right)
	\left( \sum_{\theta_2\in\faces{k}{\tau}} \sgn{\tau}{\theta_2} \phi(\theta_2)\right)
%%%%%
\\&=	\sum_{\tau\in \faces{k+1}{X}} \sum_{\sigma\in\faces{k}{\tau}} \phi(\sigma)^2 \\&+\sum_{\tau\in \faces{k+1}{X}} \sum_{\theta_1\in\faces{k}{\tau}}
	\sum_{\substack{\theta_2\in\faces{k}{\tau},\\ \theta_2\neq \theta_1}} \sgn{\tau}{\theta_1}  \sgn{\tau}{\theta_2} \phi(\theta_1)\phi(\theta_2)	
%%%
\\&=	\sum_{\sigma\in \faces{k}{X}} \deg(\sigma) \phi(\sigma)^2 \\&+\sum_{\tau\in \faces{k+1}{X}} \sum_{\theta_1\in\faces{k}{\tau}}
\sum_{\substack{\theta_2\in\faces{k}{\tau},\\ \theta_2\neq \theta_1}} \sgn{\tau}{\theta_1}  \sgn{\tau}{\theta_2} \phi(\theta_1)\phi(\theta_2).
\end{align*}
%For $\tau\in \faces{k+1}{X}$ and $\theta_1, \theta_2 \in \faces{k}{\tau}$ such that $\theta_1\neq \theta_2$, we have  $|\theta_1\cap\theta_2|=k$. Take $\eta\in \faces{k-1}{X}$ such that $\theta_1\cap\theta_2=\eta$, and take $v,w\in V$ such that $\theta_1=v\eta$ and $\theta_2=w\eta$ (all the equalities are as unordered sets). Then we have $\tau=vw\eta$, therefore $vw\in\lk(X,\eta)$.
%On the other hand, given $\eta\in\faces{k}{X}$ and $vw\in\lk(X,\eta)$ we can recover $\theta_1=v\eta$, $\theta_2=w\eta$ and $\tau=vw\eta$. So we have a bijection
%Also,
%\[
%	\sum_{\tau\in \faces{k+1}{X}} \sum_{\sigma\in \faces{k}{\tau}}\phi(\sigma)^2
%	=\sum_{\sigma\in \faces{k}{X}} \degree{X}(\sigma)\phi(\sigma)^2,
%\]
%\improvement[inline]{Explain this better? or just don't explain...}
Now look at the map
\[
\left\{ (\eta,v,w) : \substack{\eta\in \faces{k-1}{X},\\
	v,w\in V,\,\,\, v\neq w,\\ \, vw\in\lk(X,\eta)}\right\}
\to
\left\{ (\tau,\theta_1,\theta_2) : \substack{\tau\in \faces{k+1}{X},\\
	\theta_1,\theta_2\in\faces{k}{\tau}, \, \theta_1\neq \theta_2}\right\}
\]
defined by $(\eta,v,w)\mapsto (vw\eta,v\eta,w\eta)$. For each $(\tau,\theta_1,\theta_2)$ in the codomain, let $\eta=\theta_1\cap\theta_2$, $\{v\}=\theta_1\setminus\theta_2$ and $\{w\}=\theta_2\setminus\theta_1$. $(\eta,v,w)$ is the unique element sent to $(\tau,\theta_1,\theta_2)$. So the map is a bijection, therefore we obtain
\begin{align*}
	\norm{\cobound{k}{\phi}}^2
	&=\sum_{\sigma\in \faces{k}{X}} \degree{X}(\sigma) \phi(\sigma)^2 \\
%%%%
	& + \sum_{\eta\in X(k-1)}\sum_{v\in V} \sum_{\substack{w\in V\setminus\{v\},\\vw \in\lk(X,\eta)}} \sgn{vw\eta}{v\eta}\sgn{vw\eta}{w\eta}\phi(v\eta)\phi(w\eta)\\
%%%%%%%%5
	&=\sum_{\sigma\in X(k)} \degree{X}(\sigma) \phi(\sigma)^2- 2 \sum_{\eta\in X(k-1)} \sum_{vw\in \lk(X,\eta)}\phi(v\eta)\phi(w\eta).	
\end{align*}
%\todo[inline]{When I say $vw\in\lk(X,\eta)$ this implies that $v\neq w$. Is this clear?}
\end{proof}

\begin{claim}
	\label{claim:dlocal}
	For $\phi\in \cochains{k}{X}$
	\begin{align*}
	\sum_{u\in V} \norm{\cobound{k-1}{\phi_u}}^2 &=
	\sum_{\sigma\in \faces{k}{X}} \sum_{\tau\in \faces{k-1}{\sigma}}\degree{X}(\tau) \phi(\sigma)^2\\
%%%%%%%%%%%%%
	&-2k \sum_{\tau\in \faces{k-1}{X}} \sum_{vw\in \lk(X,\tau)} \phi(v\tau)\phi(w\tau)\\
%%%%%%%%%%%%%%
	&-2 \sum_{\eta\in \faces{k-2}{X}}\sum_{vw\in \lk(X,\eta)} \sum_{\substack{u\in\lk(X,v\eta)\cap\lk(X,w\eta)\\u\notin \lk(X,vw\eta)}}\phi(vu\eta)\phi(wu\eta).
	\end{align*}
\end{claim}
\begin{proof}
First we apply Claim \ref{claim:dk} to $\phi_u \in\cochains{k-1}{X}$:
\[
	\norm{\cobound{k-1}{\phi_u}}^2= \sum_{\tau\in \faces{k-1}{X}} \degree{X}(\tau) \phi_u(\tau)^2
	-2 \sum_{\eta\in \faces{k-2}{X}} \sum_{vw\in \lk(X,\eta)} \phi_u(v\eta)\phi_u(w\eta).
\]
Summing over all vertices we obtain
\begin{align*}
	\sum_{u\in V} \normsquare{d_{k-1}\phi_u}
	&=\sum_{u\in V} \sum_{\tau\in X(k-1)} \degree{X}(\tau) \phi_u(\tau)^2\\
	&-2 \sum_{u \in V} \sum_{\eta\in X(k-2)} \sum_{vw\in \lk(X,\eta)} \phi_u(v\eta)\phi_u(w\eta)\\
%%%%
	&=\sum_{u\in V} \sum_{\substack{\tau\in\faces{k-1}{X}\cap\lk(X,u)}} \degree{X}(\tau) \phi(u\tau)^2\\
%%%%%%%%%
	&-2 \sum_{\eta\in X(k-2)} \sum_{vw\in \lk(X,\eta)}\sum_{u\in \lk(X,v\eta)\cap\lk(X,w\eta)} \phi(vu\eta)\phi(wu\eta)\\
%%%%%%
   &=\sum_{\sigma\in X(k)} \sum_{\tau\in \sigma(k-1)} \degree{X}(\tau) \phi(\sigma)^2\\
%%%%%%
	&-2 \sum_{\eta\in X(k-2)} \sum_{vw\in \lk(X,\eta)}\sum_{u\in \lk(X,v\eta)\cap\lk(X,w\eta)} \phi(vu\eta)\phi(wu\eta).
\end{align*}
%We have a one-to-one correspondence between pairs
% \[
%\left\{ (u,\tau) : \substack{u\in V,\\
%\tau\in\faces{k-1}{X}\cap\lk(X,u)}\right\}
%\leftrightarrow
%\left\{ (\sigma,\tau) : \substack{\sigma\in\faces{k}{X},\\
%\tau\in\faces{k-1}{\sigma}}\right\},
%\]
%(defined by sending $(u,\tau)$ to  $(\sigma,\tau)$ where $\sigma=u\tau$).
Let $\eta\in\faces{k-2}{X}$, $vw\in\lk(X,\eta)$, and $u\in \lk(X,v\eta)\cap\lk(X,w\eta)$. We split into two different cases: $u\in \lk(X,vw\eta)$ or $u\notin\lk(X,vw\eta)$.	
Assume $u\in \lk(X,vw\eta)$, and let  $\tau=u\eta$. Then we have $vw\in \lk(X,\tau)$. This defines a map
\[
	\left\{ (\eta,vw,u) : \substack{\eta\in \faces{k-2}{X},\\
	vw\in\lk(X,\eta), u\in \lk(X,vw\eta)}\right\}
	\to
	\left\{ (\tau,vw) : \substack{\tau\in \faces{k-1}{X},\\
		vw\in\lk(X,\tau)}\right\}.
\]
Each pair $(\tau,vw)$ has a preimage of size $k$ (these are the tuples $(\tau\setminus{u},vw,u)$ for each $u\in\tau$). Therefore we obtain
\begin{align*}
	\sum_{u\in V}&\norm{\cobound{k-1}{\phi_u}}^2
%%%%
	=\sum_{\sigma\in X(k)} \sum_{\tau\in \sigma(k-1)} \degree{X}(\tau) \phi(\sigma)^2\\
	&-2 \sum_{\eta\in X(k-2)} \sum_{vw\in \lk(X,\eta)}\sum_{u\in\lk(X,vw\eta)} \phi(vu\eta)\phi(wu\eta)\\
	&-2 \sum_{\eta\in X(k-2)} \sum_{vw\in \lk(X,\eta)}\sum_{\substack{u\in \lk(X,v\eta)\cap\lk(X,w\eta)\\ u\notin \lk(X,vw\eta)}} \phi(vu\eta)\phi(wu\eta)\\
%%%%%%
	&=\sum_{\sigma\in X(k)} \sum_{\tau\in \sigma(k-1)} \degree{X}(\tau) \phi(\sigma)^2
	-2k \sum_{\tau\in X(k-1)} \sum_{vw\in \lk(X,\tau)} \phi(v\tau)\phi(w\tau)\\
	&-2 \sum_{\eta\in X(k-2)} \sum_{vw\in \lk(X,\eta)}\sum_{\substack{u\in \lk(X,v\eta)\cap\lk(X,w\eta)\\ u\notin \lk(X,vw\eta)}} \phi(vu\eta)\phi(wu\eta).
	\end{align*}
\end{proof}
\begin{remark} If $X$ is a clique complex and $u\in \lk(X,v\eta)\cap\lk(X,w\eta)$ for $\eta\in\faces{k-2}{X}$ and $vw\in\lk(X,\eta)$, then all the $1$-dimensional faces of the simplex $uvw\eta$ belong to $X$, therefore $uvw\eta\in X$ (i.e. $u\in\lk(X,vw\eta)$). Therefore in this case the last term of the previous equation vanishes (see \cite[Claim 3.2]{aharoni2005eigenvalues}).
%\improvement[inline]{Explicar mejor y poner reference a la mishbaa...}
\end{remark}

\begin{claim}[{see \cite[Claim 3.3]{aharoni2005eigenvalues}}]
\label{claim:dkstar}
For $\phi\in C^k(X)$
\[
\sum_{u\in V} \norm{\bound{k-2}{\phi_u}}^2= k \norm{\bound{k-1}{\phi}}^2.
\]
\end{claim}
\begin{proof}
%By Equation \eqref{eq:boundary}:
\[
	\norm{d_{k-1}^*\phi}^2=
%%%%%%
	\sum_{\tau\in X(k-1)} \bound{k-1}{\phi}(\tau)^2=
%%%%%%
	\sum_{\tau\in X(k-1)}\left(\sum_{v\in\lk(X,\tau)} \phi(v\tau)\right)^2.
\]
%Thus plugging $\phi_u$ %in equation (\ref{eq:dkstar}) 
%	and summing over all vertices we obtain
Similarly,
\begin{multline*}
	\sum_{u\in V} \norm{\bound{k-2}\phi_u}^2
%%%%%%	
	=\sum_{u\in V} \sum_{\eta\in \faces{k-2}{X}} \left(\sum_{v\in \lk(X,\eta)} \phi_u(v\eta)\right)^2\\
%%%%%%
	=\sum_{\eta\in X(k-2)}\sum_{u\in\lk(X,\eta)} \left(\sum_{v\in \lk(X,u\eta)} \phi(uv\eta)\right)^2\\
%%%%%%
	=k\sum_{\tau\in X(k-1)} \left(\sum_{v\in \lk(X,\tau)} \phi(v\tau)\right)^2
%%%%%
	=k \norm{\bound{k-1}\phi}^2.
\end{multline*}
\end{proof}

Let $\mata{k}:C^k(X)\to C^k(X)$ be the linear transformation whose matrix representation in the standard basis is
\[
\matrixrepel{\mata{k}}{\sigma}{\tau}=
\begin{cases}
(\mis{\sigma\cup\tau}-2)\cdot\sgn{\sigma}{\sigma\cap\tau}\cdot\sgn{\tau}{\sigma\cap\tau}
& \mbox{if } \substack{|\sigma\cap\tau|=k,\\ \sigma\cup\tau\notin X(k+1)},\\
0 & \mbox{otherwise.}
\end{cases}
\]

\begin{claim}
\label{claim:Amatrix}
For $\phi\in C^k(X)$
\[
	\inner{\mata{k}\phi}{\phi}=
	 2 \sum_{\eta\in \faces{k-2}{X}}\sum_{vw\in \lk(X,\eta)} \sum_{\substack{u\in\lk(X,v\eta)\cap\lk(X,w\eta),\\	u\notin\lk(X,vw\eta)}}\phi(vu\eta)\phi(wu\eta).
\]
\end{claim}

\begin{proof}
\begin{multline*}
	\inner{\mata{k}\phi}{\phi}\\
%%%%%
	=\sum_{\tau\in \faces{k}{X}} \sum_{\substack{\sigma\in \faces{k}{X},\\ |\sigma\cap\tau|=k,\\ \sigma\cup\tau\notin X(k+1)}}  (\mis{\sigma\cup\tau}-2)\sgn{\sigma}{\sigma\cap\tau}\cdot\sgn{\tau}{\sigma\cap\tau} \phi(\tau) \phi(\sigma)\\
%%%%%
	=\sum_{\theta\in \faces{k-1}{X}} \sum_{v\in \lk(X,\theta)} \sum_{\substack{w\in \lk(X,\theta),\\ vw\theta\notin \faces{k+1}{X}}} (\mis{vw\theta}-2)\sgn{v\theta}{\theta}\cdot\sgn{w\theta}{\theta} \phi(v\theta)\phi(w\theta)\\
%%%%%
	=\sum_{\theta\in \faces{k-1}{X}} \sum_{v\in \lk(X,\theta)} \sum_{\substack{w\in \lk(X,\theta),\\ vw\theta\notin \faces{k+1}{X}}} (\mis{vw\theta}-2)\phi(v\theta)\phi(w\theta).
%=\sum_{\eta\in X(k-2)} \sum_{v\in\lk(X,\eta)} \sum_{w\in \lk(X,\eta)} \sum_{\substack{u\in \lk(X,v\eta)\cap\lk(X,w\eta)\\ T(uvw\eta)=\{uvw\}}} \phi(vu\eta)\phi(wu\eta).
\end{multline*}
Let $<$ be an order on the vertices of $X$. Look at the map
\[
	\left\{ (\eta,u,vw) : \substack{\eta\in \faces{k-2}{X},\,
		vw\in\lk(X,\eta),\\ u\in\lk(X,v\eta)\cap\lk(X,w\eta),\\
		u\notin\lk(X,vw\eta)}\right\}
	\to
	\left\{ (\theta,v,w) : \substack{\theta\in \faces{k-1}{X},\\
		v,w\in\lk(X,\theta), \, v<w,\\ vw\theta\notin\faces{k+1}{X}}\right\}
\]
defined by  $(\eta,u,vw)\mapsto(u\eta,v,w)$.
Note that for any $(\eta,u,vw)$ in the domain, we must have $u,v,w\in\mi{uvw\eta}$. 
Let $(\theta,v,w)$ in the codomain, and let $u\in \mi{vw\theta}\setminus\{v,w\}$ and $\eta=\theta\setminus\{u\}$. Then $vw\in\lk(X,\eta)$ (since $vw\eta$ doesn't contain $u$, therefore can't contain any missing face). Similarly,
$u\in\lk(X,v\eta)\cap\lk(X,w\eta)$, but $u\notin \lk(X,vw\eta)$ (otherwise $uvw\eta=vw\theta\in \faces{k+1}{X}$). Therefore $(\eta,u,vw)$ is in the preimage of $(\theta,v,w)$. Hence $(\theta,v,w)$ has a preimage of size $\mis{vw\theta}-2$.
So we have
\begin{multline*}
	\inner{\mata{k}\phi}{\phi}=\sum_{\theta\in X(k-1)} \sum_{v\in \lk(X,\theta)} \sum_{\substack{w\in \lk(X,\theta),\\ vw\theta\notin X(k+1)}} (\mis{vw\theta}-2)\phi(v\theta)\phi(w\theta)\\
%%%%%
=2\sum_{\theta\in X(k-1)} \sum_{v\in \lk(X,\theta)} \sum_{\substack{w\in \lk(X,\theta),\\ v<w,\\ vw\theta\notin X(k+1)}} (\mis{vw\theta}-2)\phi(v\theta)\phi(w\theta)\\
%%%%%
=2\sum_{\eta\in X(k-2)}  \sum_{vw \in \lk(X,\eta)} \sum_{\substack{u\in \lk(X,v\eta)\cap\lk(X,w\eta),\\ u\notin\lk(X,vw\eta)}}\phi(vu\eta)\phi(wu\eta).
\end{multline*}
	
\end{proof}

\begin{proof}[Proof of Proposition \ref{prop:laplacian_identity}]
	
Let $\phi\in\cochains{k}{X}$. By Claim \ref{claim:Amatrix} we have
\begin{multline*}
%\label{eq:Bmatrix}
\inner{\matb{k}\phi}{\phi}=\sum_{\sigma\in \faces{k}{X}} \left(k\degree{X}(\sigma)-\sum_{\tau\in \faces{k-1}{\sigma}} \degree{X}(\tau)\right)\phi(\sigma)^2\\
+ 2 \sum_{\eta\in \faces{k-2}{X}}\sum_{vw\in \lk(X,\eta)} \sum_{\substack{u\in\lk(X,v\eta)\cap\lk(X,w\eta)\\u\notin\lk(X,vw\eta)}}\phi(vu\eta)\phi(wu\eta).
\end{multline*}		
By Claims \ref{claim:dk} and \ref{claim:dlocal} we obtain
\begin{multline*}
%\label{eq:dkvsdklocal}
k\norm{\cobound{k}{\phi}}^2\\=
\sum_{u\in V}\norm{\cobound{k-1}{\phi_u}}^2
+\sum_{\sigma\in \faces{k}{X}}\left( k\degree{X}(\sigma)- \sum_{\tau\in \faces{k-1}{\sigma}}\degree{X}(\tau)\right) \phi(\sigma)^2 \\+2 \sum_{\eta\in \faces{k-2}{X}}\sum_{vw\in \lk(X,\eta)} \sum_{\substack{u\in\lk(X,v\eta)\cap\lk(X,w\eta)\\u\notin \lk(X,vw\eta)}}\phi(vu\eta)\phi(wu\eta)\\
= \sum_{u\in V}\norm{\cobound{k-1}{\phi_u}}^2+\inner{\matb{k}\phi}{\phi}.
\end{multline*}
Then by the previous equation and Claim \ref{claim:dkstar}, we obtain
\begin{align*}
k\inner{\lap{k}\phi}{\phi}&=
%%%%%
k\inner{\bound{k}\cobound{k}\phi+\cobound{k-1}\bound{k-1}\phi}{\phi}
%%%%
=k \norm{\cobound{k}\phi}^2+k \norm{\bound{k-1}\phi}^2\\
%%%%
&=\sum_{u\in V}\norm{\cobound{k-1}\phi_u}^2
+\inner{\matb{k}\phi}{\phi}
+\sum_{u\in V} \norm{\bound{k-2}\phi_u}^2\\
%%%%%
&=\sum_{u\in V} \inner{\lap{k-1}\phi_u}{\phi_u}
+\inner{\matb{k}\phi}{\phi}.
\end{align*}
Substracting $(d-1)\inner{\lap{k}\phi}{\phi}$ from both sides of the equation we get
\begin{multline*}
(k-d+1)\inner{\lap{k}\phi}{\phi}= \sum_{u\in V} \inner{L_{k-1}\phi_u}{\phi_u}-\inner{((d-1)\lap{k}-\matb{k})\phi}{\phi}
\\
=\sum_{u\in V} \inner{L_{k-1}\phi_u}{\phi_u}-\inner{\matc{k}\phi}{\phi}.
\end{multline*}
\end{proof}

%Now we prove Proposition \ref{prop:eigenbound}. 
%By the definition of $\matb{k}$ and Claim \ref{claim:lapmatrix} we have
For the proof of Proposition \ref{prop:eigenbound} we will need the next result, which follows from the definition of $\matb{k}$ and Claim \ref{claim:lapmatrix}.
\begin{claim}
\label{claim:C_k}
The matrix representation of $\matc{k}$ in the standard basis %$\standardbasis{X(k)}$ 
is
\[
	\matrixrepel{\matc{k}}{\sigma}{\tau}=
	\begin{cases}
	\substack{\sum_{\eta\in\faces{k-1}{\sigma}}\degree{X}(\eta) -(k-d+1) \degree{X}(\sigma) +(d-1)(k+1)} & \mbox{if } \sigma=\tau,\\
	(d+1-\mis{\sigma\cup\tau})\sgn{\sigma}{\sigma\cap\tau}\cdot\sgn{\tau}{\sigma\cap\tau} & \mbox{if } \substack{|\sigma\cap\tau|=k,\\ \sigma\cup\tau\notin X(k+1),}\\
	0 & \mbox{otherwise.}
	\end{cases}
\]
\end{claim}
%\begin{proof}
%	Follows by the definition of $B_k$, Claim \ref{claim:lapmatrix} and by the observation we made earlier, that for any $\sigma,\tau\in X(k)$ such that $|\sigma\cap\tau|=k$ and $\sigma\cup\tau\notin X(k+1)$, then $2\leq\mis{\sigma\cup\tau}\leq d+1$.
%\end{proof}

\begin{proof}[Proof of Proposition \ref{prop:eigenbound}]
Let $\miscomplex{r}$ be the $(k+1)$-dimensional simplicial complex on vertex set $V$, with full $k$-skeleton, whose $(k+1)$-dimensional faces are the simplices $\eta \in \binom{V}{k+2}$ such that $\mis{\eta}=r$. By Claim \ref{claim:lapmatrix}, we have
\begin{multline}
	\label{eq:r_laplacian}
	\matrixrepel{\lappm{k}{\miscomplex{r}}{+}}{\sigma}{\tau}=
	\begin{cases}
	\degree{\miscomplex{r}}(\sigma) & \mbox{if } \sigma = \tau, \\
	-\sgn{\sigma}{\sigma \cap \tau}\cdot\sgn{\tau}{\sigma \cap \tau}
	& \mbox{if } \left|\sigma \cap \tau \right|=k,  \mis{\sigma \cup \tau}=r, \\
	0 & \mbox{otherwise.}
	\end{cases}	
	\\
	=\begin{cases}
	\left| \set{v\in V :\, \mis{v\sigma}=r} \right| & \mbox{if } \sigma = \tau, \\
	-\sgn{\sigma}{\sigma \cap \tau}\cdot\sgn{\tau}{\sigma \cap \tau}
	& \mbox{if } \left|\sigma \cap \tau \right|=k,  \mis{\sigma \cup \tau}=r, \\
	0 & \mbox{otherwise.}
	\end{cases}
\end{multline}	
Denote by $\tempmat{k}{r}$ the principal submatrix of $\matrixrep{\lappm{k}{\miscomplex{r}}{+}}$ obtained by keeping only the rows and columns corresponding to simplices in $X(k)$. $\tempmat{k}{r}$ is a positive semidefinite matrix (as a principal submatrix of a positive semidefinite matrix).
	
Define a new matrix
\[
M_k= \matrixrep{\matc{k}}+ \sum_{r=2}^{d} (d+1-r)\tempmat{k}{r}.
\]
For a matrix $A$, denote by $\lambda_{\text{max}}(A)$ the largest eigenvalue of $A$. 
Since $\tempmat{k}{r}$ is positive semidefinite it follows that $\lambda_{\text{max}}(-\tempmat{k}{r})\leq 0$ for all $2\leq r\leq d$ and therefore
\begin{equation}
\label{eq:ck_vs_mk}
\eigc{k}=\lambda_{\max}(\matc{k})
\leq \lambda_{\text{max}}(M_k)+\sum_{r=2}^{d} (d+1-r) \lambda_{\text{max}}(-\tempmat{k}{r})
\leq \lambda_{\text{max}}(M_k).
\end{equation}	
By equation \eqref{eq:r_laplacian}, Lemma \ref{lem:missing_faces} and  Claim \ref{claim:C_k} we see that the matrix $M_k$ is diagonal, and

\begin{multline*}
	%\label{eq:D_k}
	(M_k)_{\sigma,\sigma}=
	\sum_{\eta\in\faces{k-1}{\sigma}}\degree{X}(\eta) -(k-d+1) \degree{X}(\sigma) +(d-1)(k+1)
	\\	+\sum_{r=2}^{d} (d+1-r)\cdot \left| \set{v\in V :\, \mis{v\sigma}=r} \right|.
\end{multline*}
	%Therefore it suffices to show that $D_k(\sigma,\sigma)\leq p\cdot n$ for any $\sigma\in X(k)$.	
Let $\sigma\in X(k)$. We can write
\[
	\degree{X}(\sigma)= |\set{v\in V: v\in\lk(X,\sigma)}|
\]
and
\[
	k+1=|\set{v\in V:\, v\in\sigma}|,
\]
and by Lemma \ref{lemma:countdegrees}
\begin{multline*}
	\sum_{\eta\in\faces{k-1}{\sigma}} \degree{X}(\eta)
	=\left|\set{v\in V: v\in \sigma}\right|+(k+1)\cdot\left|\set{v\in V: v\in\lk(X,\sigma)}\right| \\
	+\sum_{r=2}^{d+1} (r-1)\cdot\left|\set{v\in V: \mis{v\sigma}=r}\right|.
\end{multline*}
Hence,
\begin{multline*}
	(M_k)_{\sigma,\sigma}=
	d\cdot\left|\set{v\in V: v\in \sigma}\right|+
	d\cdot\left|\set{v\in V: v\in \lk(X,\sigma)}\right|\\
	+\sum_{r=2}^{d+1} d\cdot\left|\set{v\in V: \mis{v\sigma}=r}\right|
	\leq d\cdot |V|= d n.
\end{multline*}
Therefore $\lambda_{\text{max}}(M_k)\leq d n$, so by inequality \eqref{eq:ck_vs_mk}: $\eigc{k}\leq d n$.
\end{proof}

\section{Vector domination}
\label{sec:vector_domination}

In this section we study the vector domination number $\Gamma(X)$ of a simplicial complex $X$, leading up to the proof of Theorem \ref{thm:connectivity_bound} that provides an upper bound on $\Gamma(X)$ in terms of the homological connectivity of $X$. First we prove Proposition \ref{prop:gamma_vs_gamma}, relating $\Gamma(X)$ to the total domination number $\tilde{\gamma}(X)$.

\begin{proof}[Proof of Proposition \ref{prop:gamma_vs_gamma}]
	Let $S$ be a totally dominating set in $X$. Let $\sigma\in\setsforrep{X}=\binom{V}{d-1}$. Let $f_{\sigma}$ be the characteristic vector of $S\setminus\sigma$. Define  $\alpha_{\sigma}=\frac{1}{d} f_{\sigma}$ if $\sigma\subset S$, and  $\alpha_{\sigma}=0$ otherwise. Then for every vector representation $P$ of $X$ and every $w\in V$ we have 
	\[
	\sum_{\sigma\in\setsforrep{X}} \sum_{v\in V} \alpha_{\sigma}(v) \vrep{\sigma}(v)\cdot \vrep{\sigma}(w)= \sum_{\sigma \in \binom{S}{d-1}} \sum_{v\in S\setminus\sigma} \frac{1}{d} \vrep{\sigma}(v)\cdot \vrep{\sigma}(w).
	\]
	$S$ is totally dominating, therefore there is some $\tau\subset S$ such that $\tau\in X$ but $w\tau\notin X$. Since all the missing faces are of dimension $d$ we must have $|\tau|\geq d$, and by taking a subset if necessary we may assume $|\tau|=d$. For every $\sigma\in\binom{\tau}{d-1}$, let $u$ be the unique vertex in $\tau\setminus\sigma$. Then $wu\sigma=w\tau$ is a missing face of $X$, thus $\vrep{\sigma}(u)\cdot\vrep{\sigma}(w)\geq 1$. Hence
	\[
	\sum_{\sigma \in \binom{S}{d-1}} \sum_{v\in S\setminus\sigma} \frac{1}{d} \vrep{\sigma}(v)\cdot \vrep{\sigma}(w) \geq   \sum_{\sigma \in \binom{\tau}{d-1}} \frac{1}{d}= d\cdot\frac{1}{d}=1.
	\]
	So $\sum_{\sigma\in\setsforrep{X}} \alpha_{\sigma}\vrep{\sigma}\vrep{\sigma}^T\geq 1$, therefore $\{\alpha_{\sigma}\}_{\sigma\in \setsforrep{X}}$ is dominating for $P$. So we have
	\[
	|P|\leq \sum_{\sigma\in\setsforrep{X}} \sum_{v\in V} \alpha_{\sigma}(v)=
	\sum_{\sigma\in \binom{S}{d-1}}\sum_{v\in S\setminus \sigma} \frac{1}{d}=
	\binom{|S|}{d-1} \frac{|S|-d+1}{d}= \binom{|S|}{d}.
	\]  
	Therefore $\Gamma(X)\leq \binom{\tot{X}}{d}$.
\end{proof}

Let $X$ be a simplicial complex. For each $i\in\missingdims{X}$, let $X_i$ be the complex whose missing faces are $\missfaces{i}{X}$. Note that $X_i$ has full $(i-1)$-dimensional skeleton and $X=\cap_{i\in \missingdims{X}} X_i$.

We want to bound the spectral gaps of $X$ by the spectral gaps of the complexes $X_i$. We will need the following lemma:
\begin{lemma}
	\label{claim:intersection_eigen}
	Let $A_1,\ldots,A_m$ be simplicial complexes on vertex set $V$, where $|V|=n$. Then
	\[
	\mineig{k}{\cap_{i=1}^m A_i}\geq \sum_{i=1}^m \mineig{k}{A_i}-(m-1)n.
	\]
\end{lemma}
\begin{proof}
We argue by induction on $m$. For $m=1$ the statement is trivial. Assume $m=2$. For any complex $C$ on vertex set $V$ containing $A_1\cap A_2$, denote by $\tilde{L}_k(C)$ the principal submatrix of $[\laplacian{k}{C}]$ obtained by keeping only the rows and columns corresponding to simplices of $A_1\cap A_2$. 

Let $\lambda_{\min}(\tilde{L}_k(C))$ and $\lambda_{\max}(\tilde{L}_k(C))$ be respectively the minimal and maximal eigenvalues of $\tilde{L}_k(C)$.
	
We have $\lambda_{\min}(\tilde{L}_k(C))\geq \mineig{k}{C}$ and $\lambda_{\max}(\tilde{L}_k(C))\leq \lambda_{\max}(\laplacian{k}{C})\leq n$ (by Lemma \ref{lemma:max_eig}).

It is easy to check by Claim \ref{claim:lapmatrix} that
\[
	[\laplacian{k}{A_1\cap A_2}]= \tilde{L}_k(A_1)+\tilde{L}_k(A_2)-\tilde{L}_k(A_1\cup A_2).
\]
Therefore
\begin{multline*}
	\mineig{k}{A_1\cap A_2} \geq \lambda_{\min}(\tilde{L}_k(A_1))
	+\lambda_{\min}(\tilde{L}_k(A_2))-\lambda_{\max}(\tilde{L}_k(A_1\cup A_2))\\
	\geq \mineig{k}{A_1}+\mineig{k}{A_2}-n.
\end{multline*}
For $m>2$ we get by the case $m=2$ and the induction hypothesis
\begin{multline*}
	\mineig{k}{\cap_{i=1}^m A_i}\geq \mineig{k}{A_1}+\mineig{k}{\cap_{i=2}^m A_i} -n
	\\\geq \mineig{k}{A_1} +\left( \sum_{i=2}^m \mineig{k}{A_i} -(m-2)n\right) -n
	=\sum_{i=1}^m \mineig{k}{A_i} -(m-1)n.
\end{multline*}
\end{proof}

For $i\in \missingdims{X}$ let $Y_i$ be the $i$-dimensional complex on vertex set $V$ with full $(i-1)$-dimensional skeleton whose $i$-dimensional faces are the sets in $\missfaces{i}{X}$. Denote the maximal eigenvalue of $\lappm{i-1}{Y_i}{+}$ by $\maxeig{i}{X}$.

\begin{claim}
	\label{claim:miu_vs_lambdamax}
	For all $i\in \missingdims{X}$ 
	\[
	\mineig{i-1}{X_i}= n- \maxeig{i}{X}.
	\]
\end{claim}
\begin{proof}
By Claim \ref{claim:lapmatrix} we have
\begin{align*}
\matrixrepel{\lap{i-1}^{+}(Y_i)}{\sigma}{\tau}&=
\begin{cases}
\degree{Y_i}(\sigma) & \mbox{if } \sigma = \tau,\\
-\sgn{\sigma}{\sigma \cap \tau}\cdot\sgn{\tau}{\sigma \cap \tau}
& \mbox{if } \substack{\left|\sigma \cap \tau \right|=i-1, \\ \sigma \cup \tau \in \faces{i}{Y_i},} \\
0 & \mbox{otherwise.}
\end{cases}\\
%%%%%%%%%%%%%%%%%%%
&=\begin{cases}
n-i-\degree{X_i}(\sigma) & \mbox{if } \sigma = \tau,\\
-\sgn{\sigma}{\sigma \cap \tau}\cdot\sgn{\tau}{\sigma \cap \tau}
& \mbox{if } \substack{\left|\sigma \cap \tau \right|=i-1, \\ \sigma \cup \tau \notin \faces{i}{X_i},} \\
0 & \mbox{otherwise.} 
\end{cases}\\
&=\begin{cases}
n-\matrixrepel{\laplacian{i-1}{X_i}}{\sigma}{\tau} & \mbox{if } \sigma = \tau,\\
-\matrixrepel{\laplacian{i-1}{X_i}}{\sigma}{\tau} & \mbox{otherwise.}\\
\end{cases}
\end{align*}
Therefore  $\lappm{i-1}{Y_i}{+}=n\identity- \laplacian{i-1}{X_i}$. So every eigenvector of $\laplacian{i-1}{X_i}$ with eigenvalue $\lambda$ is an eigenvector of of $\lappm{i-1}{Y_i}{+}$ with eigenvalue $n-\lambda$. In particular, $n-\mineig{i-1}{X_i}$ is the largest eigenvalue of $\lappm{i-1}{Y_i}{+}$.
\end{proof}

\begin{claim}
\label{claim:miu_vs_lambdamax2}
For $k\geq 0$,
\[
	\mineig{k}{X}
		\geq n- \sum_{i\in\missingdims{X}} \binom{k+1}{i} \maxeig{i}{X}.
\]
\end{claim}
\begin{proof}
By Lemma \ref{claim:intersection_eigen} we obtain
\[
\mineig{k}{X}= \mineig{k}{\cap_{i\in \missingdims{X}} X_i}
\geq \sum_{i\in \missingdims{X}} \mineig{k}{X_i} - (|\missingdims{X}|-1)n.
\]
Applying Claim \ref{claim:induction} to each of the complexes $X_i$ (note that $\maxmis(X_i)=i$ and $X_i$ has full $(i-1)$-dimensional skeleton) we get
\[
\mineig{k}{X}\geq \sum_{i\in \missingdims{X}} \left[\binom{k+1}{i} \mineig{i-1}{X_i}- \left(\binom{k+1}{i}-1\right)n\right] - (|\missingdims{X}|-1)n.
\]
Then by Claim \ref{claim:miu_vs_lambdamax}
\begin{multline*}
%\label{eq:mineig_vs_maxeigs}
\mineig{k}{X}\geq \sum_{i\in \missingdims{X}} \left[\binom{k+1}{i}(n-\maxeig{i}{X})- \left(\binom{k+1}{i}-1\right)n\right] - (|\missingdims{X}|-1)n\\
=n-\sum_{i\in \missingdims{X}} \binom{k+1}{i}\maxeig{i}{X}.
\end{multline*}
\end{proof}

\begin{claim}
\label{claim:eigenhom2}
\[
	\sum_{i\in \missingdims{X}} \binom{\eta(X)}{i}\maxeig{i}{X} \geq n.
\]
\end{claim}
\begin{proof}
Let $k$ be the integer such that
\[
	\sum_{i\in \missingdims{X}} \binom{k-1}{i} \maxeig{i}{X} < n \leq	\sum_{i\in \missingdims{X}} \binom{k}{i}\maxeig{i}{X}.
\]
Let $j\leq k-2$. By Claim \ref{claim:miu_vs_lambdamax2},
\[
\mineig{j}{X}\geq n-\sum_{i\in \missingdims{X}} \binom{j+1}{i}\maxeig{i}{X} > 0,
\]
therefore by Corollary \ref{cor:hodge} we have $\homology{j}{X}=0$. So $\eta(X)\geq k$, thus
	\[
	\sum_{i\in \missingdims{X}} \binom{\eta(X)}{i}\maxeig{i}{X} \geq
	\sum_{i\in \missingdims{X}} \binom{k}{i}\maxeig{i}{X} \geq n.
	\]
\end{proof}

\begin{claim}
\label{claim:pluslapnorm}
Let $i\in\missingdims{X}$. Then for $\phi\in\cochains{i-1}{Y_i}$,
\[
	\inner{\lappm{i-1}{Y_i}{+}\phi}{\phi}\leq \sum_{\sigma\in\binom{V}{i-1}}\sum_{vw\in\lk(Y_i,\sigma)} (\phi(v\sigma)-\phi(w\sigma))^2.
\]
\end{claim}
\begin{proof}
\begin{multline*}
\sum_{\sigma\in\faces{i-2}{Y_i}} \sum_{vw\in\lk(Y_i,\sigma)} (\phi(v\sigma)-\phi(w\sigma))^2\\
%%%%%
=\sum_{\sigma\in \faces{i-2}{Y_i}} \sum_{ v\in \lk(Y_i,\sigma)} \degree{Y_i}(v\sigma) \phi(v\sigma)^2 -2 \sum_{\sigma\in \faces{i-2}{Y_i}} \sum_{vw\in \lk(Y_i,\sigma)} \phi(v\sigma)\phi(w\sigma)\\
%%%%%
=i\cdot\sum_{\eta\in \faces{i-1}{Y_i}} \degree{Y_i}(\eta) \phi(\eta)^2 -2 \sum_{\sigma\in \faces{i-2}{Y_i}} \sum_{vw\in \lk(Y_i,\sigma)} \phi(v\sigma)\phi(w\sigma).
\end{multline*}
By Claim \ref{claim:dk} we have
\begin{multline*}
\inner{\lappm{i-1}{Y_i}{+}\phi}{\phi}
%%%%%%
=\norm{d_{i-1}\phi}^2 \\
%%%%%
=\sum_{\eta \in \faces{i-1}{Y_i}} \degree{Y_i}(\eta) \phi(\eta)^2-
2\sum_{\sigma\in \faces{i-2}{Y_i}} \sum_{vw\in \lk(Y_i,\sigma)} \phi(v\sigma)\phi(w\sigma).
\end{multline*}
Hence
\begin{multline*}
\inner{\lappm{i-1}{Y_i}{+}\phi}{\phi}
%%%%%%
=\sum_{\sigma\in\faces{i-2}{Y_i}} \sum_{vw\in\lk(Y_i,\sigma)} (\phi(v\sigma)-\phi(w\sigma))^2 \\
-(i-1)\cdot \sum_{\eta \in \faces{i-1}{Y_i}} \degree{Y_i}(\eta) \phi(\eta)^2
%%%%%%
\leq \sum_{\sigma\in\faces{i-2}{Y_i}} \sum_{vw\in\lk(Y_i,\sigma)} (\phi(v\sigma)-\phi(w\sigma))^2.
\end{multline*}
$Y_i$ has full $(i-1)$-dimensional skeleton, therefore $\faces{i-2}{Y_i}=\binom{V}{i-1}$. Thus
\[
\inner{\lappm{i-1}{Y_i}{+}\phi}{\phi}\leq \sum_{\sigma\in\binom{V}{i-1}} \sum_{vw\in\lk(Y_i,\sigma)} (\phi(v\sigma)-\phi(w\sigma))^2
\]
\end{proof}

\begin{claim}
\label{claim:eigenrep}
Let $P$ be a vector representation of $X$. Then for all $i\in \missingdims{X}$
\[
	\maxeig{i}{X}\leq i \cdot \max_{\substack{\sigma\in \binom{V}{i-1},\, v\in V}} \left( P_{\sigma}(v)\cdot \sum_{w\in V} P_{\sigma}(w)\right).
\]
\end{claim}

\begin{proof}
Let $\phi\in \cochains{i-1}{Y_i}$. For $\sigma\in\faces{i-2}{Y_i}=\binom{V}{i-1}$ and $v,w\in V\setminus\sigma$, $v\neq w$, we have, by the definition of $P$, $\vrep{\sigma}(v)\cdot\vrep{\sigma}(w)\geq 1$ if $vw\in\lk(Y_i,\sigma)$, and $\vrep{\sigma}(v)\cdot\vrep{\sigma}(w)\geq 0$ otherwise. Therefore we obtain
\begin{align}
\label{eq:vecrep1}
&\sum_{\sigma\in\binom{V}{i-1}} \sum_{vw\in\lk(Y_i,\sigma)} (\phi(v\sigma)-\phi(w\sigma))^2 \nonumber\\
%%%%%%
&\leq \frac{1}{2}\sum_{\sigma\in\binom{V}{i-1}} \sum_{v,w\in V\setminus\sigma} (\phi(v\sigma)-\phi(w\sigma))^2  \vrep{\sigma}(v)\cdot\vrep{\sigma}(w) \nonumber\\
%%%%%
&=\sum_{\sigma\in \binom{V}{i-1}} \sum_{v\in V\setminus\sigma} \phi(v\sigma)^2 \vrep{\sigma}(v)\cdot \sum_{w\in V\setminus\sigma} \vrep{\sigma}(w) 
-\sum_{\sigma\in\binom{V}{i-1}} \norm{\sum_{v\in V\setminus \sigma} \phi(v\sigma) \vrep{\sigma}(v)}^2 \nonumber\\
%%%%%
&\leq \sum_{\sigma\in\binom{V}{i-1}} \sum_{v\in V\setminus\sigma} \phi(v\sigma)^2 \vrep{\sigma}(v)\cdot \sum_{w\in V\setminus\sigma} \vrep{\sigma}(w)   \nonumber\\
%%%%%%%
&\leq 
\left(\sum_{\sigma\in \binom{V}{i-1}} \sum_{v\in V\setminus\sigma} \phi(v\sigma)^2 \right)\cdot \max_{\substack{\sigma\in \binom{V}{i-1},\, v\in V\setminus \sigma}} \vrep{\sigma}(v)\cdot \sum_{w\in V\setminus\sigma} \vrep{\sigma}(w).
\end{align}
Since $Y_i$ has full $(i-1)$-dimensional skeleton, we have
\begin{multline}
\label{eq:vecrep2}
\sum_{\sigma\in\binom{V}{i-1}} \sum_{v\in V\setminus\sigma} \phi(v\sigma)^2
%%%%%%%%%%
=\sum_{\sigma\in\faces{i-2}{Y_i}} \sum_{v\in\lk(Y_i,\sigma)} \phi(v\sigma)^2\\
%%%%%%%%%%%%%
= i \sum_{\eta\in \faces{i-1}{Y_i}} \phi(\eta)^2
%%%%%%
=  i \norm{\phi}^2.
\end{multline}
Combining \eqref{eq:vecrep1},\eqref{eq:vecrep2} and Claim \ref{claim:pluslapnorm} we obtain
\begin{multline*}
\inner{\lappm{i-1}{Y_i}{+}\phi}{\phi}\leq
\sum_{\sigma\in\binom{V}{i-1}} \sum_{vw\in\lk(Y_i,\sigma)} (\phi(v\sigma)-\phi(w\sigma))^2\\
%%%%%%%%%%
\leq i \norm{\phi}^2 
\cdot \max_{\substack{\sigma\in \binom{V}{i-1},\, v\in V\setminus \sigma}} \vrep{\sigma}(v)\cdot \sum_{w\in V\setminus\sigma} \vrep{\sigma}(w) \\
%%%%%%
\leq  i \norm{\phi}^2\cdot \max_{\sigma\in \binom{V}{i-1},\, v\in V} \vrep{\sigma}(v)\cdot \sum_{w\in V} \vrep{\sigma}(w).
\end{multline*}
Thus
\[
\maxeig{i}{X}= \max_{0\neq\phi\in \cochains{i-1}{Y_i}} \frac{\inner{L_{i-1}^+\phi}{\phi}}{\|\phi\|^2} \leq i\cdot \max_{\sigma\in \binom{V}{i-1} ,v\in V} \left(P_{\sigma}(v)\cdot \sum_{w\in V} P_{\sigma}(w)\right).
\]
\end{proof}

\begin{lemma}
\label{lem:linear_duality}
Let $P$ be a vector representation of $X$. Then
\[
|P|= \max\set{\alpha\cdot\textbf{1}:\quad \alpha\geq 0, \, \alpha \vrep{\sigma} \vrep{\sigma}^T\leq \textbf{1} \quad \forall \sigma\in \setsforrep{X} }.
\]
\end{lemma}
\begin{proof}
Let $\sigma_1,\ldots,\sigma_m$ be all the sets in $\setsforrep{X}$. For each $i\in \set{1,2,\ldots,m}$ let $A_i=\vrep{\sigma_i}\vrep{\sigma_i}^T\in \Rea^{|V|\times|V|}$. Note that $A_i=A_i^T$. % for all $i\in\set{1,2,\ldots,m}$.
Define the matrix
\[	
	A=\left( A_1 | A_2 | \cdots | A_m \right)^T\in \Rea^{(m|V|)\times |V|}.
\]

Let $x\in \Rea^{m|V|}$. Write $x=\left(\alpha_{\sigma_1}| \alpha_{\sigma_2}| \cdots | \alpha_{\sigma_m}\right)$, where $\alpha_{\sigma_i}\in \Rea^{|V|}$ for each $i\in \set{1,2,\ldots,m}$. We have
\[
	xA=\sum_{i=1}^m \alpha_{\sigma_i} A_i = \sum_{\sigma\in\setsforrep{X}}\alpha_{\sigma}\vrep{\sigma}\vrep{\sigma}^T,
\]
therefore
\begin{align*}
|P|&=\min\set{\sum_{\alpha\in\setsforrep{X}} \alpha_{\sigma} \cdot \textbf{1} : \, \alpha_{\sigma}\geq 0 \,\, \forall \sigma\in \setsforrep{X},
\sum_{\sigma\in\setsforrep{X}}\alpha_{\sigma}\vrep{\sigma}\vrep{\sigma}^T \geq \textbf{1} } \\
%%%%%%%%%%%%%%%%%%
&=\min\set{ x\cdot\textbf{1} : \, x\geq 0, xA\geq \textbf{1}}.
\end{align*}
By linear programming duality
\[
|P|=\max\set{ y\cdot\textbf{1} : y\geq 0, yA^T\leq \textbf{1}}.
\]
But $yA^T=\left(yA_1| yA_2 | \cdots | yA_m\right)$, so $yA^T\leq \textbf{1}$ if and only if $y A_i\leq \textbf{1}$ for all $i\in\set{1,2,\ldots,m}$. Therefore
\[
	|P|=\max\set{y\cdot\textbf{1} : y\geq 0, \, y\vrep{\sigma}\vrep{\sigma}^T\leq \textbf{1}\quad \forall \sigma\in\setsforrep{X}}.
\]	
\end{proof}

Let $\ZZ_{+}$ denote the positive integers, and $\QQ_{+}$ the positive rationals. Let $a\in\ZZ_{+}^V$ and
\[
	\multiplier{a}{V}=\set{ (v,i) : \, v\in V, 1\leq i\leq a(v)}.
\]
Define the projection $\pi: \multiplier{a}{V}\to V$ by $\pi((v,i))=v$, and let 
\[
	\multiplier{a}{X}=\pi^{-1}(X)=\set{\sigma\subset \multiplier{a}{V} : \, \pi(\sigma)\in X}.
\]
The missing faces of $\multiplier{a}{X}$ are the sets $\sigma\subset \multiplier{a}{V}$ such that $|\pi(\sigma)|=|\sigma|$ and $\pi(\sigma)$ is a missing face of $X$.

$\pi$ induces an homotopy equivalence between $\multiplier{a}{X}$ and $X$ (see \cite[Lemma 2.6]{lovasz1996topological}), therefore $\conn{\multiplier{a}{X}}=\conn{X}$.

\begin{proof}[Proof of Theorem \ref{thm:connectivity_bound}]
Let $P=\set{\vrep{\sigma}}_{\sigma\in \setsforrep{X}}$ be a vector representation of $X$. Let $\alpha\in \QQ_{+}^V$ such that $\alpha \vrep{\sigma} \vrep{\sigma}^T\leq \textbf{1}$ for all $\sigma\in \setsforrep{X}$. Write $\alpha=a/k$ where $k\in \ZZ_{+}$ and $a\in \ZZ_{+}^V$.
Denote $N=|\multiplier{a}{V}|= \sum_{v\in V} a(v)$.
For $\sigma\in\setsforrep{\multiplier{a}{X}}$ and $(v,j)\in\multiplier{a}{V}$ define 
%$Q_{\sigma}((v,j))=\vrep{\pi(\sigma)}(v)$ if $|\pi(\sigma)|=|\sigma|$ and $Q_{\sigma}((v,j))=0$ otherwise.
\[
	Q_{\sigma}((v,j))=\begin{cases}
	\vrep{\pi(\sigma)}(v) & \text{ if } |\pi(\sigma)|=|\sigma|,\\
	0 & \text{ otherwise.}
	\end{cases}
\]
$Q=\set{Q_{\sigma} : \, \sigma\in \setsforrep{\multiplier{a}{X}}}$ is a vector representation of $\multiplier{a}{X}$: Let $\sigma\in\setsforrep{\multiplier{a}{X}}$ of size $r-1$, and let $\tilde{v}=(v,i),\tilde{u}=(u,j)\in\multiplier{a}{V}$ such that $\tilde{u}\tilde{v}\sigma\in\missfaces{r}{\multiplier{a}{X}}$. Then $\pi(\tilde{u}\tilde{v}\sigma)=uv\pi(\sigma)\in \missfaces{r}{X}$. In particular $|\pi(\sigma)|=|\sigma|$, therefore, since $P$ is a representation of $X$,
\[
	Q_{\sigma}(\tilde{v})\cdot Q_{\sigma}(\tilde{u})= \vrep{\pi(\sigma)}(v)\cdot\vrep{\pi(\sigma)}(u)\geq 1.
\]	
Let $r\in\missingdims{X}$. By Claim \ref{claim:eigenrep}
\begin{align*}
\maxeig{r}{\multiplier{a}{X}} &\leq  r \cdot \max_{\sigma\in\binom{\multiplier{a}{V}}{r-1}, (v,j)\in\multiplier{a}{V}}
\left( Q_{\sigma}((v,j))\cdot \sum_{(w,k)\in \multiplier{a}{V}} Q_{\sigma}((w,k))\right)\\
%%%%%%%%%%%%%%
&=r \cdot \max_{\tau\in\binom{V}{r-1},v\in V}\left( \vrep{\tau}(v)\cdot \sum_{w\in V} a(w) \vrep{\tau}(w)\right)
	\leq r\cdot k.
\end{align*}
By Claim \ref{claim:eigenhom2} we obtain
\[
\sum_{r\in\missingdims{\multiplier{a}{X}}} \binom{\conn{\multiplier{a}{X}}}{r} r\cdot k\geq
%%%%%%%%%%%%%%%%%%
\sum_{r\in\missingdims{\multiplier{a}{X}}} \binom{\conn{\multiplier{a}{X}}}{r} \maxeig{r}{\multiplier{a}{X}}  
\geq N.
\]
Therefore
\[
	\alpha\cdot\textbf{1} = \frac{1}{k} \sum_{v\in V} a(v)= \frac{N}{k}
	\leq \sum_{r\in\missingdims{\multiplier{a}{X}}} r \binom{\eta(\multiplier{a}{X})}{r}=\sum_{r\in\missingdims{X}} r \binom{\eta(X)}{r}.
\]
Thus by Lemma \ref{lem:linear_duality}
\begin{align*}
|P|&= \max\set{\alpha\cdot\textbf{1}:\, \alpha\geq 0, \, \alpha \vrep{\sigma} \vrep{\sigma}^T\leq \textbf{1} \, \forall \sigma\in \setsforrep{X} }\\
%%%%%%%%%%%
&=\sup \set{\alpha\cdot\textbf{1}:\, \alpha\in\QQ_{+}^V, \, \alpha \vrep{\sigma} \vrep{\sigma}^T\leq \textbf{1} \, \forall \sigma\in \setsforrep{X} }\leq
\sum_{r\in\missingdims{X}} r \binom{\conn{X}}{r},
\end{align*}
therefore $\Gamma(X)\leq \sum_{r\in\missingdims{X}} r \binom{\conn{X}}{r}$.
\end{proof}

For the proof of Theorem \ref{thm:generalhalltype} we need the following Hall-type condition for the existence of colorful simplices, which appears in \cite{aharoni2000hall,meshulam2001clique}, and more explicitly in \cite{meshulam2003domination}:
\begin{proposition}
\label{prop:colorfulsimplex}
Let Z be a simplicial complex on vertex set $W=\cupdot_{i=1}^m W_i$. If for all $\emptyset\neq I\subset \set{1,2,\ldots m}$
\[
	\eta(Z[\cupdot_{i\in I} W_i])\geq |I|
\]
then $Z$ contains a colorful simplex.
\end{proposition}

\begin{proof}[Proof of Theorem \ref{thm:generalhalltype}]
Let $\emptyset\neq I \subset \set{1,2,\ldots,m}$. By Theorem $\ref{thm:connectivity_bound}$ we have
\[
	\sum_{r\in\missingdims{X[\cupdot_{i\in I}V_i]}} r\binom{\eta(X[\cupdot_{i\in I}V_i]))}{r}\geq\Gamma(X[\cup_{i\in I}V_i])> \sum_{r\in\missingdims{X[\cupdot_{i\in I}V_i]}} r \binom{|I|-1}{r},
\]
therefore
	\[
	\eta(X[\cupdot_{i\in I}V_i]))>|I|-1.
	\]
Thus by Proposition $\ref{prop:colorfulsimplex}$ $X$ has a colorful simplex.
\end{proof}

\section{Colorful sets in general position}
\label{sec:matroids}
Let $M$ be a matroid of rank $d+1$ on vertex set $V$. Let $\tilde{M}$ be the simplicial complex on vertex set $V$ whose simplices are the subsets $S\subset V$ in general position with respect to $M$.
The missing faces of $\tilde{M}$ are the dependent sets $S\subset V$ with $|S|\leq d+1$ such that any $|S|-1$ points in $S$ are independent in $M$.
%\todo[inline]{Add explanation about missing faces? The missing faces are the subsets $S\subset V$ of size $k$ such that $\matroidrank(S)=k-1$ and any subset of size $k-1$ has closure equal to $\matroidclosure(S)$ (i.e. any subset of size $k-1$ has rank $k-1$, i.e. any such subset is independent).}

%So the missing faces of $\tilde{M}$ are the subsets $S\subset V$ of size $k$ that are contained in a flat of $M$ of rank $k-1$ but aren't contained in any flat of lower rank (that is, any subset of $S$ of size $k-1$ is an independent set), for $2\leq k\leq d+1$.-THIS IS WRONG

\begin{claim}
\label{claim:varphistar_vs_gamma}
For $U\subset V$,
\[
	\varph^*(U)\leq d\cdot \Gamma(\tilde{M}[U]).
\]
\end{claim}
\begin{proof}
We construct a vector representation of the complex $\tilde{M[U]}$. Let $1\leq r\leq d$ and let $\mathcal{F}_r$ be the set of flats of $M$ of rank $r$.
	
Let $\sigma\in\setsforrep{\tilde{M}[U]}$ with $|\sigma|=r-1$, and let $v\in U$. Define $\vrep{\sigma}(v)\in \Rea^{\mathcal{F}_r}$ by
\[
	\vrep{\sigma}(v)(F)=\begin{cases}
	1 & \text{ if } \matroidclosure(v\sigma)=F,\\
	0 & \text{ otherwise.}
	\end{cases}
\]
For $v,w\in U$, if $vw\sigma$ is a missing face of $\tilde{M}[U]$ of dimension $r$ then $vw\sigma$ lies in a flat of rank $r$, which is spanned by any $r$ points in $vw\sigma$. In particular $\matroidclosure(v\sigma)=\matroidclosure(w\sigma)\in \mathcal{F}_r$, therefore
\[
	\vrep{\sigma}(v)\cdot \vrep{\sigma}(w)=1.
\]
Hence $P$ is a vector representation of $\tilde{M}[U]$.

Let $f: U\to \Rea_{\geq 0}$ be a function in fractional general position with $\sum_{v\in U} f(v)=\varph^*(U)$. Define $\alpha\in\Rea^{U}$ by $\alpha(v)=f(v)/d$. 

Let $w\in U$, and let $F=\matroidclosure(w\sigma)$. If $F \notin \mathcal{F}_{r}$ then $\vrep{\sigma}(w)=0$, therefore $\sum_{v\in U} \alpha(v)\vrep{\sigma}(v)\cdot\vrep{\sigma}(w)=0\leq 1$. If $F\in\mathcal{F}_r$ then
 
\[
	\sum_{v\in U} \alpha(v) \vrep{\sigma}(v)\cdot \vrep{\sigma}(w)
%%%%%%%%
=\sum_{\substack{v\in U,\\ \matroidclosure(v\sigma)=F}} \alpha(v)
%%%%%%%%%%%%
=\frac{1}{d} \sum_{\substack{v\in U,\\ \matroidclosure(v\sigma)=F}} f(v)\leq 1.
\]
So $\alpha \vrep{\sigma}\vrep{\sigma}^T \leq \textbf{1}$ for each $\sigma\in\setsforrep{\tilde{M}[U]}$, therefore by Lemma \ref{lem:linear_duality}
\[
	\Gamma(\tilde{M}[U])\geq |P|\geq \alpha\cdot \textbf{1}= \frac{\varph^*(U)}{d}.
\]

\end{proof}

\begin{proof}[Proof of Theorem \ref{thm:myHMSstar}]
Let $\emptyset\neq I\subset \set{1,2,\ldots,m}$.
%Note that the induced subcomplex $X[\cupdot_{i\in I} V_i]$ is just $X(\cupdot_{i\in I} V_i)$, the complex on vertex set $\cupdot_{i\in I} V_i$ whose simplices are the subsets in general position.
By Claim \ref{claim:varphistar_vs_gamma}
\[
\Gamma(\tilde{M}[\cupdot_{i\in I} V_i])\geq \frac{\varph^*(\cupdot_{i\in I} V_i)}{d}> \sum_{r=1}^d r\binom{|I|-1}{r}.
\]
Thus by Theorem \ref{thm:generalhalltype} there is a colorful simplex of $\tilde{M}$, i.e. a colorful subset of $V$ in general position.
\end{proof}

\begin{proof}[Proof of Theorem \ref{thm:myHMS}]
Let $\emptyset\neq I\subset \set{1,2,\ldots,m}$.
Assume $|I|\leq d+1$.
The $d$-dimensional skeleton of $\tilde{M}[\cupdot_{i\in I} V_i]$ is $M[\cupdot_{i\in I}V_i]$, therefore for all $0\leq k\leq d-1$
\[
	\cohomology{k}{\tilde{M}[\cupdot_{i\in I} V_i]}=\cohomology{k}{M[\cupdot_{i\in I} V_i]}.
\]
$M$ is a matroid, therefore $\cohomology{k}{M[\cupdot_{i\in I} V_i]}=0$ for $0\leq k\leq \matroidrank(\cupdot_{i\in I} V_i)-2$ (see \cite{bjorner1995topological}).
So $\conn{\tilde{M}[\cupdot_{i\in I} V_i]}\geq \matroidrank(\cupdot_{i\in I} V_i)$.
But
\[
\matroidrank(\cupdot_{i\in I} V_i)=\min\set{d+1,\varph(\cupdot_{i\in I} V_i)},\] so if $\varph(\cupdot_{i\in I} V_i)>|I|-1$, then $\conn{\tilde{M}[\cupdot_{i\in I} V_i]}> |I|-1$.

Assume now that $|I|\geq d+2$. If $\varph(\cupdot_{i\in I} V_i)>d \sum_{r=1}^d r\binom{|I|-1}{r}$, then, by inequality \eqref{eq:fractionalphi}, $\varph^*(\cupdot_{i\in I} V_i)>d \sum_{r=1}^d r\binom{|I|-1}{r}$, and therefore by Theorem \ref{thm:connectivity_bound} and Claim \ref{claim:varphistar_vs_gamma}
\[
\sum_{r=1}^d r \binom{\conn{\tilde{M}[\cupdot_{i\in I}V_i]}}{r} \geq \Gamma(\tilde{M}[\cupdot_{i\in I} V_i])\geq \frac{\varph^*(\cupdot_{i\in I} V_i)}{d}> \sum_{r=1}^d r\binom{|I|-1}{r},
\]
so $\conn{\tilde{M}[\cupdot_{i\in I}V_i]}>|I|-1$.
Therefore by Proposition \ref{prop:colorfulsimplex} there is a colorful subset of $V$ in general position.
\end{proof}

\section*{Acknowledgment}
This paper was written as part of my M. Sc. thesis, under the supervision of Professor Roy Meshulam. I thank Professor Meshulam for his guidance, and for his helpful comments and suggestions.

\bibliographystyle{plain}
\bibliography{biblio}

\end{document}